\documentclass[12pt,a4paper,reqno]{amsart}
\usepackage{amsmath}
\usepackage{amsfonts}
\usepackage{amssymb}

\newtheorem{theorem}{Theorem}[section]
\newtheorem{corollary}{Corollary}[section]
\newtheorem{lemma}{Lemma}[section]
\newtheorem{proposition}{Proposition}[section]
\newtheorem{definition}{Definition}[section]
\newtheorem{conjecture}{Conjecture}[section]
\newtheorem{assumption}{Assumption}[section]

\def\A{{\mathbb A}}
\def\C{{\mathbb C}}

\def\E{{\mathcal E}}
\def\F{{\mathbb F}}
\def\G{{\mathbb G}}
\def\H{{\mathbb H}}

\def\K{{\mathcal K}}

\def\O{{\mathcal O}}

\def\Q{{\mathbb Q}}

\def\cS{{\mathcal S}}

\def\Z{{\mathbb Z}}

\DeclareMathOperator{\Gal}{Gal}
\DeclareMathOperator{\res}{res}
\DeclareMathOperator{\ext}{ext}

\DeclareMathOperator{\sign}{sign}
\DeclareMathOperator{\Val}{Val}
\DeclareMathOperator{\Hom}{Hom}
\DeclareMathOperator{\Reg}{Reg}
\DeclareMathOperator{\Pic}{Pic}
\DeclareMathOperator{\id}{id}
\DeclareMathOperator{\Ver}{Ver}
\begin{document}
\title{Generalized Stark formulae over function fields}
\author{Ki-Seng Tan}
\address{Department of Mathematics\\
National Taiwan University\\
Taipei 10764, Taiwan}
\email{tan@math.ntu.edu.tw}
\thanks{The author was supported in part by the National
Science Council of Taiwan (R.O.C.), NSC91-2115-M-002-001,
NSC93-2115-M-002-007. }

 \maketitle
\begin{abstract}
We establish formulae of Stark type for the Stickelberger elements
in the function field setting. Our result generalizes a work of Hayes
and a conjecture of Gross. It is used to deduce a $p$-adic version of
Rubin-Stark Conjecture and Burns Conjecture.
\end{abstract}

\begin{section}{Introduction}\label{s:intro}
In this paper, we study Stickelberger elements related to abelian
extensions over global function fields. Our main result is Theorem \ref{thm:stark1},
which generalizes a theorem of Hayes (\cite{h}). In a way the theorem puts together
conjectures of Gross, Rubin and Stark. And we will show that it implies a $p$-adic version of
Rubin-Stark conjecture (see Theorem \ref{th:rp}
below). Furthermore,
using the theorem, we are able to deduce a $p$-adic version of a formula
conjectured by Burns (see Theorem \ref{th:burnsp} below).
%In short, we would sat that our main formula
%this formula compares the regulators obtained
%from two different Stark units.

For the purpose of having a better description of our work, we shall
review in the following paragraphs both Rubin-Stark conjecture and Burns conjecture.
But before we do so, let us first
fix some notations.

From now on, $K/k$ will be a finite abelian extension over a global
function field of characteristic $p$. We assume that the extension
is unramified outside a given finite set $S$ of places of $k$.
And we fix another finite non-empty set $T$ of places of $k$ such that
$T\cap S=\emptyset$. The notation $K'$ will be used to denote a subfield of $K$
containing $k$. Also, $S(K')$ (resp. $T(K')$) will denote the set of places
of $K'$ sitting over $S$ (resp. $T$).
Let $\F_q$ be the constant field of $k$, and put
$\Gamma=\Gal(K/k)$, $\Gamma'=\Gal(K'/k)$.

The analytical side of Rubin-Stark conjecture involves the equivariant $L$-function
which interpolates $L$-functions at each character. Recall that for each $\chi\in
\hat \Gamma$ the modified Artin $L$-function over $k$ is defined as
(\cite{g1})
\begin{equation}\label{eq:ldef}
L_{S,T}(\chi,s)=\prod_{v\in T} (1-\chi([v])\cdot
N(v)^{1-s})\prod_{v\not\in S} (1-\chi([v])\cdot N(v)^{-s}
)^{-1}.
\end{equation}
Here $[v]$ is the Frobenius element at $v$ and
$N(v)=q^{\deg(v)}$ is the norm of $v$. Because in a global field places of the same degree
always form a finite set, one can easily deduce
that the above infinite product is expanded in a unique way as a formal
power series in $q^{-s}$. It is well-known that this formal power
series is in fact a polynomial in $q^{-s}$ (\cite{t1}).
Applying the theory of Fourier transforms, we see
that there is a polynomial
$\Theta_{\Gamma,S,T}(s)\in \C[\Gamma][q^{-s}]$
such that for every $\chi\in {\hat \Gamma}$
\begin{equation}\label{e:Theta}
\chi(\Theta_{\Gamma,S,T}(s))=L_{S,T}(\chi,s).
\end{equation}
This $\Theta_{\Gamma,S,T}(s)$ is called the modified equivariant
$L$-function.
From (\ref{eq:ldef}), we are able to express
$\Theta_{\Gamma,S,T}$ as an infinite product. Namely,
\begin{equation}\label{e:Thetaprod}
\Theta_{\Gamma,S,T}(s)=\prod_{v\in T} (1-[v]\cdot
q^{\deg(v)(1-s)})\prod_{v\not\in S} (1-[v]\cdot q^{-\deg(v)s})^{-1},
\end{equation}
In particular, this implies that $\Theta_{\Gamma,S,T}$ is actually
an element in $\Z[\Gamma][q^{-s}]$.

Now let us start to describe the arithmetic side of Rubin-Stark conjecture. This
will involve various regulator maps related to units groups.
Consider $\O_{S(K')}$, the ring of $S(K')$-integers of $K'$, and let $\O_{S(K')}^*$ be its
units group.
\begin{definition}\label{d:stu}
Define $U(K')$ to be the kernel of the reduction modulo $T(K')$
$$\O_{S(K')}^*\longrightarrow \bigoplus_{w\in T(K')} \O_w^*/
(1+\pi_w\cdot  \O_w).$$
And define $r_{K'}=\# S(K')-1$. For simplicity, we
denote $U=U(K)$, $r=r_K$.
\end{definition}

Note that $U(K')=U^{\Gal(K/K')}$ is a free abelian group of rank $r_{K'}$
and there is an exact sequence (\cite{g1})
$$1\longrightarrow U(K')\longrightarrow \O_{S(K')}^*\longrightarrow
\prod_{w\in T(K')} \F_w^*\longrightarrow \Pic(\O_{S(K')})_{T(K')}\longrightarrow
\Pic(\O_{S(K')})\longrightarrow 1.$$
We recall that the $(S(K'),T(K'))$-class number of $K'$ is the group order
\begin{equation}\label{e:hkst}
h_{K',S(K'),T(K')}=|\Pic(\O_{S(K')})_{T(K')}|.
\end{equation}

To construct the regulator maps, we shall follow the notations and the methods used in
\cite{r}. In particular,
if $M$ is a finite $\Z[\Gamma]$-module then $\Q M$ denotes $\Q\otimes M$,
and the dual module $M^*$
is defined as
%the finite $\Z[\Gamma]$-module
$\Hom_{\Gamma} (M,\Z [\Gamma])\subset \Hom_{\Gamma} (\Q M,\Q [\Gamma])$.
Also, if $n$ is a non-negative integer, then $\Lambda^n M$ denotes the
$n$th exterior power of $M$ in the category of $\Z
[\Gamma]$-modules.
We let
$\iota$ denote the natural map
(\cite{r}, Sec.1.2) $\iota :\Lambda^n (M^*)\longrightarrow (\Lambda^n M)^*$,
such that if
$\phi_1,...\phi_n\in M^*$ and $m_1,...,m_n\in M$, then
\begin{equation}\label{eq:iota}
\iota(\phi_1\wedge\dots \wedge \phi_n)(m_1\wedge\dots\wedge
m_n)=\det(\phi_i(m_j)).
\end{equation}
And following \cite{r}, we define
$$ \Lambda_0^n M:=\{ m\in\Q \Lambda^n M\; | \;
\iota(\phi_1\wedge\dots\wedge\phi_n)(m)\in\Z[\Gamma]
\;\text{for every}\; \phi_1,...,\phi_n\in M^*\}.$$

Now, we start to define the regulator maps.
First, let $Y(K')=\bigoplus_{w\in S(K')} \Z\cdot
w$,
and
$$X(K')=\{ \sum_{w\in S(K')} a_w\cdot w\in Y(K') | \sum_{w\in S(K')} a_w=0\}.$$
For each place $w$ of $K$,
let $\deg_w$ be the local degree map
$$\deg_w:{K}_w^*\longrightarrow K_w^*/\O_w^*\longrightarrow\Z,$$
such that if $|\; |_w$ is the normalized absolute value associated to $w$,
then
$$\log(|x|_w)=- \deg_w(x)\cdot \log(q).$$
Compose this local degree map with the
natural embedding $U\longrightarrow K_w^*$ to form
$$
\lambda_w: U \longrightarrow {K}_w^*\stackrel{\deg_w}{\longrightarrow} \Z.$$
And we define $\lambda:U(K)\longrightarrow X(K)$ to be the $\Gamma$-equivariant homomorphism
such that
$\lambda(u)=\sum_{w\in S(K)} \lambda_w(u)\cdot w$ for every
$u\in U(K)$.
Write $\lambda^{(n)}:\Lambda^n U\longrightarrow \Lambda^n X(K)$
(resp. $i^{(n)}:\Lambda^n X(K)\stackrel{i^{(n)}}{\longrightarrow} \Lambda^n Y(K)$)
for the map induced by the map $\lambda$ (resp. the inclusion
$i:X(K)\longrightarrow Y(K)$).
Then the regulator map
\begin{equation}\label{e:classreg}
\mathcal{R}_{\Psi}:\Q\Lambda^n U\longrightarrow \Q[\Gamma]
\end{equation}
associated to a $\Psi\in \Lambda^n Y(K)^*$ is defined as the one that linearly extends the map
$$\Lambda^n U\stackrel{\lambda^{(n)}}{\longrightarrow}
\Lambda^n X(K)\stackrel{i^{(n)}}{\longrightarrow}
\Lambda^n Y(K)\stackrel{\iota(\Psi)}{\longrightarrow} \Z[\Gamma].$$

It is easy to see that if $\psi_i\in Y(K)^*$, then $\mathcal{R}_{\psi_i}\in U^*$. Also,
if $\Psi=\psi_1\wedge\dots\wedge\psi_n$, then we have
\begin{equation}\label{e:etaiota}
\mathcal{R}_{\Psi}=\iota(\mathcal{R}_{\psi_1}\wedge\dots\wedge
\mathcal{R}_{\psi_n}).
\end{equation}
Consequently, we have
$ \mathcal{R}_{\Psi}(\Lambda_0^n U)\subset \Z[\Gamma]$ for every $\Psi\in \Lambda^n Y(K)^*$.
In this paper, some special elements in $Y(K)^*$ of the form $w^*$ will be used.
Recall that for each place $w$ over $K$, the element
$w^*\in Y(K)^*$ is defined (\cite{r}) such that
for every place $w'$
$$w^*(w')=\sum_{\gamma w=w'}\gamma.$$

What we called Rubin-Stark Conjecture
is the one proposed by Rubin in \cite{r} (Conjecture B$^{\prime}$), because it can be viewed
as an integral version of Stark's conjecture (\cite{s, stark1, stark2, stark3, stark4}).
In a way, the conjecture relates some derivative of the equivariant $L$-function
to certain exterior product of units arising form regular representations of $\Gamma$.
It is easy to see that the units group $U$ contains a regular representation of $\Gamma$
if and only if some place in $S$ splits completely over $K$. Thus,
for the purpose of having an interesting theory, we need to assume the following:
\begin{assumption}\label{a:a}
From now on, we assume that there exist $n$, $n\geq 1$, different places
$S_0=\{v_1,...,v_n\}\subsetneq S$
such that every place in $S_0$ splits completely over $K$,

\end{assumption}
\begin{definition}\label{d:fix}
Let $v_1,...,v_n$ be as in {\em {Assumption \ref{a:a}}}
and let $w_1,...,w_n$ be a fixed set of places of $K$ such that each $w_i$ is
sitting over the place $v_i\in S$. And let
$\eta=w_1^*\wedge\dots\wedge w_n^*$.
\end{definition}

First we note that
Assumption \ref{a:a} (together with the class number formula at $s=0$) implies
$$\Theta_{\Gamma,S,T}=a_n (q^{-s}-1)^n+\dots \in
(q^{-s}-1)\cdot \Z[\Gamma][q^{-s}].$$
And the coefficient $a_n\in\Z[\Gamma]$ will be denoted as $\Theta_{\Gamma,S,T}^{(n)}(0)$.
%Let ${F}$ be the subfield of $\C$, which is generated over $\Q$ by
%all the values of the characters of $\Gamma$, and let $\O$ be its
%ring of integers. We fix a prime ideal $\P$ of $\O$ sitting over
%$p$, and let ${F}_p$ (resp. $\O_p$) denote the $\P$-completions of ${ F}$
%(resp. of $\O$).
For each $\chi\in {\hat {\Gamma}}$ let $e_{\chi}$ be
the associated idempotent element in the group ring $\C[\Gamma]$,
and let $r_{\chi}$ denote the $\C$-dimension of the
%let $U_{\chi}$ be the
$\chi$-eigenspace of $\C\otimes_{\Z} U$.
%$F_p\otimes U$ and let $r_{\chi}=\dim_{F_p} U_{\chi}$.
Define
$$\Lambda^n_{S,T}=\{u\in \Lambda_0^n U| e_{\chi}(u)=0,\;\;
\text{for every}\;\;\chi\in {\hat \Gamma}\;\;
\text{such that}\;\;r_{\chi}>n\}.$$
Then in our
settings Rubin-Stark conjecture reads as following.
\begin{conjecture}\label{c:r}\em{(Rubin,\cite{r}, Conjecture B$^{\prime}$)}
There exists an $\epsilon\in\Lambda_{S,T}^n$ such that
\begin{equation}\label{e:rub}
\mathcal{R}_{\eta}(\epsilon)=\Theta_{\Gamma,S,T}^{(n)}(0).
\end{equation}
\end{conjecture}

We will show in Section \ref{sec:mainthm} that our main result implies
the following $p$-adic version of the conjecture. Here "$p$-adic" means tensoring things with
$\Z_{(p)}$. In particular,
$$\Z_{(p)}\Lambda_{S,T}=\Z_{(p)}\otimes \Lambda_{S,T}\subset \Q\Lambda_{S,T}.$$
\begin{theorem}\label{th:rp}
There exists an $\epsilon\in\Z_{(p)}\Lambda_{S,T}^n$ such that in $\Z_{(p)}[\Gamma]$
\begin{equation}\label{e:rp}
\mathcal{R}_{\eta}(\epsilon)=\Theta_{\Gamma,S,T}^{(n)}(0).
\end{equation}
\end{theorem}
A different proof of the theorem can be found in \cite{p05},
and a proof for the $l$-adic ($l\not= p$)
version of Rubin-Stark conjecture is given in \cite{b05}. In view of this,
over function fields, Rubin-Stark conjecture actually holds.
Now we review the conjecture of Burns.
We will follow the construction
described in \cite{b02,ha04}. The conjecture involves regulators of another type,
and we are going to define them in the follow paragraphs.

First we note that if $M$ is a $\Z[\Gamma]$-module,
then for each $\phi\in M^*$ there is a unique
$\phi^{(id)}\in\Hom_{\Z}(M,\Z)$ such that for $x\in M$
\begin{equation}\label{e:phiphiid}
\phi(x)=\sum_{\gamma\in\Gamma} \phi^{(id)}(\gamma^{-1}x)\gamma.
\end{equation}
For a place $v$ over $k$, define
\begin{equation}\label{e:lvgamma}
{\bar {\lambda}}_{v,\Gamma}:U(k)\longrightarrow k_v^*\longrightarrow
\Gamma_v \hookrightarrow\Gamma,
\end{equation}
where the first and the last arrows are natural embeddings
and the second is the norm residue map in the local class field theory.
Let $u_1,...,u_{r_k}$ be a $\Z$-basis for $U(k)$,
$v_{n+1},...,v_{r_k}$ be distinct places in $S\setminus S_0$
and $\phi_1,...,\phi_n\in U^*$. Consider the matrix $A=(a_{ij})_{1\leq i,j
\leq r_k}$ with
$$a_{ij}=\begin{cases}
\phi_i^{(id)}(u_j),\text{if} \; 1\leq i\leq n\\
{\bar {\lambda}}_{v_i,\Gamma}(u_j)-1, \text{if} \; n+1\leq i\leq r_k.
\end{cases}
$$
For each pair $i,j$ the entry $a_{ij}$ is an element in $\Z[\Gamma]$.
And it is obvious that the determinant $\det(A)$ is in  $I^{r_k-n}$ where
$I$ is the augmentation ideal of $\Z[\Gamma]$.
Up to $\pm 1$, the residue class of
$\det(A)$ modulo $I^{r_k-n+1}$ depends on neither the ordering of $v_1,...,v_{r_k}$
nor the choice of the basis $u_1,...,u_{r_k}$. We assume that the ordering of $v_1,..., v_{r_k}$
is fixed and the basis $u_1,...,u_{r_k}$ are ordered in a way such that the classical regulator
formed by them is positive. On the other hand,  the residue class of
$\det(A)$ actually depends on the exterior product
$\Phi=\iota(\phi_1\wedge\dots\wedge\phi_n)\in \iota(\Lambda^n U^*)$,
and therefore we
will denote it as $\Reg_{\Gamma}^{\Phi}$. On top of Conjecture \ref{c:r},
Burns \cite{b02} proposes the following strengthened conjecture.
%For this reason, we call it the conjecture of Rubin and Burns.
\begin{conjecture}\label{c:burns}
Assume that \em{Conjecture \ref{c:r}} holds, so that for every $\Phi\in
\iota(\Lambda^n U^*)$,
we have $\Phi(\epsilon)\in \Z[\Gamma]$. Then this element satisfies
$$\Phi(\epsilon)\equiv h_{k,S,T}\Reg_{\Gamma}^{\Phi}\pmod{I^{r_k-n+1}}.$$
\end{conjecture}
For more material related to this conjecture, see for instance
\cite{b02, b05, ha04, p99-1, p99-2, p02,r}.
The $l$-adic version (for $l\not=p$) of the conjecture
is proved in \cite{b05}, but it seems the technique used in the proof
can not be applied to cover the following $p$-adic version, which
 will be proved in Section \ref{sec:gross}.
Let $I_p$ be the augmentation ideal of $\Z_p[\Gamma]$.

\begin{theorem}\label{th:burnsp}
Let notations be as those in \em{Theorem \ref{th:rp}}.
Then for every $\Phi\in\iota(\Lambda^n U^*)$,
we have $\Phi(\epsilon)\in\Z_p[\Gamma]$
and this element satisfies
\begin{equation}\label{e:burnsp}
\Phi(\epsilon)\equiv h_{k,S,T}\Reg_{\Gamma}^{\Phi}\pmod{I_p^{r_k-n+1}}.
\end{equation}
\end{theorem}

Now we begin to describe Theorem \ref{thm:stark1}, our main result.
In short, it is a $p$-adic refinement of Theorem
\ref{th:rp}. The method for making this kind of refinement comes from \cite{g1, g2}, and the
main idea is to replace $\Z_p$ by certain Galois groups in order to construct refinements of both side of the equality
(\ref{e:rp}).
To explain it, let us start with
those degree maps $\deg_w$ which play important roles in the construction
of the regulator maps. These local degree maps together form
the global degree map $\deg:\A_K^*\longrightarrow \Z$ defined
on the ideles group $\A_K^*$.
Let $L_0=K\F_{q^{p^{\infty}}}$ be the constant
$\Z_p$-extension over $K$.
If we view $\Z_p$ as the Galois group
$\Gal(L_0/K)$
and compose the map $\deg$ with the embedding
$\Z\longrightarrow \Z_p$ which sends $1$ to the Frobenius in
$\Gal(L_0/K)$,
then we get the norm residue map $\A_K^*\longrightarrow
\Gal(L_0/K)$, and the local degree map is just the composite
$K_w^*\longrightarrow \A_K^*\longrightarrow \Gal(L_0/K)$. From this we see that
the field extension $L_0/K$ and the related norm residue maps
are implicitly used in the construction
of the previous regulator maps.

For the refinements we are going to use various Galois groups
of the form $H:=\Gal(L/K)$ where $L/K$ is a pro-$p$ abelian extension
such that $L/k$ is also abelian and unramified outside $S$ (such extension is
called admissible, see Definition \ref{d:ad}). We let
$H$ play the role of $\Z_p=\Gal(L_0/K)$
and use the related norm residue maps
to construct, for each $\Psi$, the associated refined
regulator map $\mathcal{R}_{\Psi,H}$
(Definition \ref{d:regulatormap}) which has values in the $n$th
relative augmentation
quotient associated to $H$ (Definition \ref{d:relaug}). To see that $\mathcal{R}_{\Psi,H}$
actually refines $\mathcal{R}_{\Psi}$, we only need to take $L=L_0$, because
in this situation $\mathcal{R}_{\Psi}$ can be recovered from $\mathcal{R}_{\Psi,H}$
(Lemma \ref{l:regrelation}). We would like to emphasize that $\H$, the Galois group of the maximal
admissible extension, is a direct product of countable infinite many copies of
$\Z_p$ (Lemma \ref{le:maxadmissible}). And, in a way,
Lemma \ref{le:torfree} together with
the isomorphism (\ref{eq:mapd}) says that an element in
the $n$th relative augmentation quotient associated to $\H$ can be identified as a $\Z_p[\Gamma]$-coefficient
$n$'th degree homogeneous "polynomial in countable infinite many variables". Furthermore, under this
identification, if the $\Z_p$-basis of
$\H$ is suitably arranged, then for each $\epsilon\in \Z_{(p)}\Lambda_{S,T}^{(n)}$ the value
$\mathcal{R}_{\Psi}(\epsilon)$ is just the coefficient of certain monomial in $\mathcal{R}_{\Psi,\H}(\epsilon)$.
In particular, it is fair to say that the map $\mathcal{R}_{\eta,\H}$, where $\eta$ is the one in Definition
\ref{d:fix}, carries a rich amount of information about the units group. In fact, the
universal property studied in
Section \ref{sec:unique} (see
Corollary \ref{c:projection}) tells us that most of the important information about the
integer structure of the $n$'th exterior product of the units group
can be obtained from $\mathcal{R}_{\eta,\H}$.

The refinement of the equivariant $L$-function
$\Theta_{\Gamma, S,T}$
turns out to be the Stickelberger
element $\theta_G$ (see Definition \ref{def:thetadef})
where $G=\Gal(L/k)$. For its reason, please see Lemma \ref{l:Thetatheta}.
It is somewhat a surprise since the Stickelberger element
only interpolates special values of $L$-functions while the equivariant $L$-function
interpolates the complete $L$-functions.
Lemma \ref{l:Thetatheta} also tells us that
in the case where $L=L_0$, the "$n$th derivative" $\Theta_{\Gamma, S,T}^{(n)}(0)$
can be recovered from the residue class $[\theta_{G}]_{(n,H)}$
of $\theta_G$ in the $n$th relative augmentation quotient.
That there is a unique $\epsilon$ in
$\Z_{(p)}\Lambda_{S,T}^{(n)}$ such that
$$[\theta_{G}]_{(n,H)}=\mathcal{R}_{\eta,H}(\epsilon)$$
for every admissible $H$ is exactly
the content of Theorem \ref{thm:stark1}. In the case where $H=\H$ we have an
equality between two "polynomials in infinite variables" while (\ref{e:rp}) in Theorem \ref{th:rp}
is an equality between the corresponding coefficients of certain "monomial".
And this is the reason why Theorem \ref{thm:stark1} implies Theorem \ref{th:rp}.

The proof of Theorem \ref{th:burnsp}
involves a refined class number formula proposed by Gross \cite{g1} (see Conjecture
\ref{c:gross}). What we actually
use is its $p$-adic version proved in \cite{tn} (see Theorem \ref{th:gross})
in which the congruence (\ref{e:grossp}) relates the Stickelberger element $\theta_G$
with the product of $h_{k,S,T}$ and a regulator ${\det}_G$
defined by Gross. In contrast to this, Theorem \ref{thm:stark1}
relates $\theta_G$ with the refined regulator $\mathcal{R}_{\eta,H}(\epsilon)$
which is the left-hand side of the congruence (\ref{e:burnsp}).
The main step for proving Theorem \ref{th:burnsp} is to use the aforementioned universal property
to relate the right-hand side of (\ref{e:burnsp}) to the product
$h_{k,S,T}{\det}_G$.

Finally, let us have some words about the proof of Theorem \ref{thm:stark1}.
In brief, it is based on two observations. First, we find that, via Fourier transform,
the theorem is equivalent to its twisted version, Theorem \ref{thm:stark2} in which the main
part is the congruence (\ref{e:stark2}).
And we have discovered that both the left-hand and right-hand sides of (\ref{e:stark2})
can be found as factors of the corresponding left-hand and right-hand sides of the
congruence (\ref{e:grossp}) in Theorem \ref{th:gross}. Furthermore, the two sides of
(\ref{e:grossp}) are indeed products of these kind of factors
(indexed by characters, see Proposition
\ref{prop:factorization}
and Proposition \ref{prop:thetafac}). To use (\ref{e:grossp}) to prove
(\ref{e:stark2}), we apply Fourier transforms, the universal property and the result of Hayes
for the case $n=1$.

This manuscript has evolved through several versions, since 1996.
It is a great pleasure to thank David Burns, Wen-Cheng Chi, Benedict Gross, Po-Yi Huang,
King F. Lai,
Cristian Popescu, Karl Rubin
and John Tate for stimulating discussions.
\end{section}

\begin{section}{Admissible extensions and Augmentation Quotients}
\label{sec:asaq}
In this chapter, we study admissible extensions and the properties of the
associated augmentation quotients.

\begin{subsection}{Admissible extensions}\label{sec:ad}

\begin{definition}\label{d:ad}
An abelian extension $L/K$ and its Galois group $H=\Gal(L/K)$ are admissible if
the followings are satisfied:
\begin{enumerate}
\item The extension $L/k$ is abelian and is unramified outside $S$.
\item The extension $L/K$ is a pro-$p$ extension.
\end{enumerate}
\end{definition}

Throughout this paper, we will fix an admissible extension
$L/K$, and we will also fix the notations: $G=\Gal(L/k)$,
$H=\Gal(L/K)$, $\Gamma=G/H$.  Also, a subgroup of $G$ denoted as $H'$ always contains $H$, and
we always denote $K'=L^{H'}$ and $\Gamma'=G/H'$.
\end{subsection}

\begin{subsection}{The maximal admissible extension}\label{sec:maxadm}
Although there are infinitely many different admissible extensions, the theory in this paper
can be summed up to a theory for a single extension, that is, the maximal admissible extension
with respect to $K/k$ and $S$. We will denote the associated Galois group
by $\mathbb{H}$ and will first study its structure.
\begin{lemma}\label{le:maxadmissible}
The maximal admissible Galois group $\mathbb{H}$ is a direct product of
countable infinite many copies of $\Z_p$.
\end{lemma}

Before we prove the lemma, let us recall some known results
related to the local Leopoldt conjecture (see \cite{k, tn}).
\begin{lemma}\label{le:localleo}
Suppose that ${\K}$ is a global function field of characteristic
$p$ and $v$ is a place over $\K$.  If an element $u\in {\K}^*$ is
divisible by $p$ in ${\K}_v^*$, then it is divisible by $p$ in
${\K}^*$.
\end{lemma}
As a consequence, we have the following .
\begin{lemma}\label{le:leop} {\em{(\cite{k})}} Suppose that ${\K}$ is a
global function field of characteristic $p$ and ${\cS}$ is a
finite set of places of ${\K}$. Then the Galois group of the
maximal pro-$p$ abelian extension over ${\K}$ unramified outside
${\cS}$ is a direct product of countable infinite many copies of
$\Z_p$.
\end{lemma}
\begin{proof} (of Lemma \ref{le:maxadmissible})

Let $\Gamma=\Gamma_p\oplus\Gamma_0 $ be the natural decomposition of $\Gamma$
into the $p$-part, $\Gamma_p$, and the non-$p$-part, $\Gamma_0$. Suppose that
${\G}$ is the Galois group over ${k}$ of the maximal pro-$p$
abelian extension unramified outside ${S}$. Then ${\G}$ is an extension of
$\Gamma_p$ which is viewed as a quotient group of $\Gamma$. Let
${\mathcal{H}}={\ker} (\G \longrightarrow \Gamma_p)$
be the kernel of the natural quotient map. Then ${\mathcal{H}}$ is
isomorphic to $\mathbb{H}$.
\par
By Lemma \ref{le:leop}, ${\G}$ is a direct product of countable
infinite many copies of $\Z_p$, and so is ${\H}$.
\end{proof}
\end{subsection}

\begin{subsection}{Group rings and augmentation
ideals}\label{sec:gring} For the rest of this chapter, we will study
group rings with
various coefficient rings together with two types of augmentation ideals and the associated augmentation quotients.

Let $R$ be an integral domain finite over
$\Z$ or $\Z_p$. If $C$ is the fraction field of $R$ and $M$ is an
$R$-module, then we use $C M$ to denote $C\otimes_R M$.
\begin{definition}\label{d:gring} For a pro-finite group
$\mathcal{H}$, let $R[\mathcal{H}]$ be the projective limit of
$ R[\Delta]$, where $\Delta$ runs through all the finite quotient groups
of $\mathcal{H}$. Also, for every positive integer $n$,
let $I_R(\mathcal{H})^n$ be the projective limit of $ I_R(\Delta)^n,$
where $I_R(\Delta)^n$ is the $n$th power of the augmentation ideal
$I_R(\Delta)$. We call respectively $I_R(\mathcal{H})^n$ and
$I_R(\mathcal{H})^n/I_R(\mathcal{H})^{n+1}$ the $n$th augmentation ideal
and the $n$th augmentation quotient of $R[\mathcal{H}]$.
\end{definition}

%Then $R[H']$ is nothing but the ring of $R$-valued measures on the
%topological group $H'$ and $I_R(H')$ is the ideal of measures on
%$H'$ with total measure zero.
For the rest of the paper, if $\Xi :\mathcal{H}_1\longrightarrow \mathcal{H}_2$
is a group homomorphism, then we will also use $\Xi$ to denote the
induced homomorphisms on the group rings and the augmentation quotients.

\begin{definition}\label{d:relaug} If $H'$ is finite, let
$I_{R,H'}$ be the kernel of the ring homomorphism $R[G]\longrightarrow R[\Gamma']$
induced from the natural quotient map $G\longrightarrow \Gamma'$.
In general, for every positive integer $n$, let $I_{R,H'}^n$ be the projective limit
of $ I_{R, H'/\mathcal{N}}^n$,
where $ \mathcal{N}$ runs through the family of all open subgroups of $H'$
contained in $H$. We call respectively $I_{R,H'}^n$ and $I_{R,H'}^n/I_{R,H'}^{n+1}$
the $n$th relative $H'$-augmentation ideal and the $n$th relative $H'$-augmentation quotient
of $R[G]$
\end{definition}

For simplicity, we let $I(H')$, $I_p(H')$, $I_{H'}$ and
$I_{p,H'}$ denote respectively $I_{\Z}(H')$, $I_{\Z_p}(H')$, $I_{\Z,H'}$ and
$I_{\Z_p, H'}$.

\begin{definition}\label{d:nthterm}
For $\xi\in I_R(H')^n\subset R[H']$, let
$[\xi]_{(n)}$ be its residue class in the augmentation quotient
$I_R(H')^n/I_R(H')^{n+1}$.
Also, for $\xi\in I_{R,H'}^n\subset R[G]$, let $[\xi]_{(n,H')}$ be its residue class  in
$I_{R,H'}^n/I_{R,H'}^{n+1}$.
\end{definition}

\begin{definition}\label{d:fring} Let $F$ be a fixed number field containing all the
values of characters of $\Gamma$ and let $\O=\O_F$ be its ring of integers.
Let $\O_p$ be the completion of $\O$ at a fixed place sitting over $p$ and let $F_p$
be its field of fractions.
Define the group rings over $F$, $F_p$ and the corresponding
augmentation ideals as follow. Let
$F[H']=F\O[H']$ and
$F_p[H']=F_p\O_p[H']$.
Also, for every positive integer $n$, let
$I_F(H')^n=F I_\O(H')^n$,
$I_{F_p}(H')^n=F_p I_{\O_p}(H')^n$,
$I_{F,H'}^n=F I_{\O,H'}^n$,
and $I_{F_p,H'}^n=F_p I_{\O_p,H'}^n$.
\end{definition}

In many situations, the structure of the augmentation quotients
can be explicitly expressed. First of all, we have the following
isomorphism (\cite{g1})
\begin{equation}\label{eq:del}
\delta_{H'}:H'\longrightarrow I(H')/I(H')^2
\end{equation}
which sends $h\in H'$ to $h-1\pmod{I(H')^2}$. Also,
if $H\simeq\Z_p^d$ for some $d$ and $R$ is either
$\Z_p$ or ${\O}_p$, then the graded ring formed by augmentation quotients
can be identified with a polynomial ring. To see this, let
$ R[[s_1,...,s_d]]$ be the ring of formal power series in $d$
variables. If ${\E}=\{\sigma_1,...,\sigma_d \}$ is a basis of $H$
over $\Z_p$ and $x_i=\sigma_i-1\in R[H]$, $i=1,...,d$, then the
map $R[H]\longrightarrow R[[s_1,...,s_d]]$, $x_i\mapsto s_i$ is
an isomorphism. Consequently, for every positive integer $n$,
\begin{equation}
 I_R(H)^n=(x_1,...,x_d)^n\simeq (s_1,...,s_d)^n,
\label{eq:ideal}
\end{equation}
and the augmentation quotient $I_R(H)^n/I_R(H)^{n+1}$ is
isomorphic to the $R$-module of $n$th degree homogeneous
polynomials in $s_1,...,s_d$. This induces an isomorphism
\begin{equation}
d_{{\E},R}:\bigoplus_{n=0}^{\infty} I_R(H)^n/I_R(H)^{n+1}
{\longrightarrow} R[s_1,...,s_d]. \label{eq:mapd}
\end{equation}
Since $F_p[H]=F_p\otimes_{{\O}_p} {\O}_p[H]$, tensoring with $F_p$, We get the induced
ring homomorphism
\begin{equation}\label{eq:mapdf}
d_{{\E},F_p}:
\bigoplus_{n=0}^{\infty} I_{F_p}(H)^n/I_{F_p}(H)^{n+1}
{\longrightarrow} F_p[s_1,...,s_d].
\end{equation}

\begin{lemma}\label{le:basechange}

\begin{enumerate}
\item[(a)]
For each non-negative
integer $n$ we have $I_p(H)^n\cap \Z[H]=I(H)^n$,
\newline
$I_{{\O}_p}(H)^n\cap {\O}[H]=I_{\O}(H)^n$
and $I_{F_p}(H)^n\cap {F}[H]=I_{F}(H)^n$.
\item[(b)] Suppose that either $H$ is finite free over $\Z_p$ or $H=\H$.
For $i=1,2$, let $A_i$ be one of the rings $\Z$, $\Z_p$, ${\O}$,
${\O}_p$, $F$ and $F_p$. If $A_1\subset A_2$, then for each $n$,
\newline
$I_{A_2}(H)^n\cap A_1[H]=I_{A_1}(H)^n$.
\item[(c)] If $H$ is finite free over $\Z_p$ or $H=\mathbb{H}$,
then the natural map
$$i:I(H)^n/I(H)^n\longrightarrow I_p(H)^n/I_p(H)^{n+1}$$ is an
isomorphism.

\end{enumerate}
\end{lemma}
\begin{proof}
\par
In Part (a), the third equality is from the second. If $H$ is
finite, then the proof of the first equality can be found in
\cite{tn}, Lemma 2.5. The second equality can be proved in a
similar way. If $H=\H$, we prove them by taking projective limits.

To prove Part (b), we first assume that $H$ is finite free over $\Z_p$.
If we are in the special case where $\{A_1,A_2\}\subset
\{\Z_p,\O_p, F_p\}$, then Part (b) is proved by using
(\ref{eq:ideal}). In general, put $A_i'=(A_i)_p$, for $i=1,2$.
Then we have
$I_{A_2'}(H)^n\cap A_1'[H]=I_{A_1'}(H)^n$.
Also, part (a) implies that
$I_{A_i'}(H)^n\cap A_i[H]=I_{A_i}(H)^n$.
These imply Part (b). The case $H=\H$ is proved by taking projective
limits.

\par
To prove Part (c), we note that by Part (b), the map $i$ is injective. First assume that $H$ is finite
free over $\Z_p$. Then by Equation (\ref{eq:ideal}),
the homomorphism $\phi$ sending
$h$ to the residue class of $h-1$ is an isomorphism from $H$ to
$I_p(H)/I_p(H)^2$. It is obvious that (see (\ref{eq:del}))
$\phi=i\circ \delta_H$. Since $\delta_H$ is an isomorphism, so is $i$.
This proves the lemma for $n=1$. For $n>1$, we observe that the
multiplication map
$$I_p/I_p^2\times I_p^{n-1}/I_p^n\longrightarrow I_p^n/I_P^{n+1}$$
maybe not surjective but its image generates the whole group (as
an abelian group). Then the surjectivity of $i$ is proved by
induction. Again, the $H=\H$ case can be proved by taking projective limits.
\end{proof}

\begin{lemma}\label{le:torfree}
For $g\in G$, let $\gamma_g\in\Gamma$
be its residue class modulo
$H$. Then for a nonnegative $n$, the homomorphism
$$
\pounds_n:  R[\Gamma]\otimes_R
I_R(H)^n/I_R(H)^{n+1}{\longrightarrow}
 I_{R,H}^n/I_{R,H}^{n+1},$$
 which sends the residue class of
$ \gamma_g\otimes (h_1-1)\cdot\dots\cdot (h_n-1)$ to that of
\newline
$g(h_1-1)\cdot\dots\cdot (h_n-1)$, $h_1,...,h_n\in H$, is an
isomorphism.
\end{lemma}
\begin{proof}
That the homomorphism $\pounds_n$ is well defined is due to the simple
fact that if $h_1,h_2\in H$, $g\in G$, then
$gh_1\equiv gh_2 \pmod{I_H}$.
The rest is obvious for the case where $H$ is finite. In general,
it is proved through projective limits.
\end{proof}
The lemma shows that the structure of the augmentation quotient
$ I_{R,H}^n/I_{R,H}^{n+1}$ actually depends only on the structures
of $\Gamma$ and $H$.

\end{subsection}

\begin{subsection}{Numerical extensions}\label{sec:numerical}
The relative $H$-augmentation quotients can be easily expressed in the following
situation.
\begin{definition}\label{d:numerical}
The extension $L/K$ and its Galois group $H$ are called numerical if
$H\simeq \Z_p$.
\end{definition}
In particular, $L/K$ is numerical, if it is
the constant $\Z_p$-extension.
If $H$ is numerical and $\sigma$ is a $\Z_p$-generator of it,
then the isomorphisms in the previous sections together form the
following isomorphism
$$\Val_{\sigma,n}=\Val_{\sigma,n,G/H}:I_H^n/I_H^{n+1}\stackrel{\pounds_n^{-1}}{\longrightarrow}
\Z[\Gamma]\otimes I(H)^n/I(H)^{n+1}
%\stackrel{id\otimes i}{\longrightarrow}
%\Z[\Gamma]\otimes I_p(H)^n/I_p(H)^{n+1}
\stackrel{d_{\sigma,\Z_p}}{\longrightarrow}
\Z_p[\Gamma].$$
Here we identify $\Z[\Gamma]\otimes I(H)^n/I(H)^{n+1}$ with
$\Z[\Gamma]\otimes I_p(H)^n/I_p(H)^{n+1}$, and we identify
an one variable homogeneous polynomial
with its coefficient. If $\sigma'=u\sigma$, $u\in\Z_p^*$, is another generator, then
$\Val_{\sigma,n}=u^n\cdot \Val_{\sigma',n}$.

If $H$ is numerical, then $G$ can be identified with $\Gamma'\times H'$
for some subgroup $H'\simeq \Z_p$ containing $H$. In this case,
$\Gamma=\Gamma'\times H'/H$.
We will relate the group ring and the $H$-augmentation quotients of the
group $G$ to those of the direct product $ \breve{G}=\Gamma\times H$.
To do so, we let $\varpi=|H'/H|$ and define
\begin{equation}\label{e:yen}
\begin{array}{rcl}
\yen: G=\Gamma'\times H' & \longrightarrow & \Gamma'\times H'/H \times H
=\Gamma\times H\\
(\gamma',h') & \mapsto & (\gamma', {\overline h'}, \varpi h')\\
\end{array}
\end{equation}

Both $G$ and $\breve{G}$ are extensions of $\Gamma$ by $H$, and
by Lemma \ref{le:torfree} we have the
associated isomorphisms
$\pounds_{n,G}:\Z[\Gamma]\otimes I(H)^n/I(H)^{n+1}\longrightarrow
I_H^n/I_H^{n+1}$ and
\newline
$\pounds_{n,\breve{G}}:\Z[\Gamma]\otimes I(H)^n/I(H)^{n+1}\longrightarrow
\breve{I}_H^n/\breve{I}_H^{n+1}$, where $\breve{I}_H^{n}$ denotes the $n$th
relative $H$-augmentation ideal of $\Z[\breve{G}]$.
The following lemma is obvious.
\begin{lemma}\label{l:yen}
Let $n$ be a nonnegative integer. An element $\xi\in\Z[G]$ is in $I_H^n$,
if and only if $\yen(\xi)$ is in $ \breve{I}_{H}^n$. In fact, if
as in {\em{Lemma} \ref{le:torfree}} we
associate the isomorphisms $\pounds_{n,G}$ and $\pounds_{n,\breve{G}}$
respectively to the groups $G$ and $\breve{G}$, then
\begin{equation}\label{e:yenpounds}
\yen\circ \pounds_{n,G}= \varpi^n\pounds_{n,\breve{G}},
\end{equation}
and
\begin{equation}\label{e:yenpounds2}
\Val_{\sigma,n,G/H}\circ \yen\circ \pounds_{n,G}=\varpi^n\Val_{\sigma,n,\breve{G}/H}\circ
\pounds_{n,\breve{G}}.
\end{equation}
\end{lemma}

\end{subsection}
\end{section}

\begin{section}{Stickelberger elements as refinements of the equivariant $L$-functions}
\label{sec:eqandst}
In this chapter, we review the definition of Stickelberger elements and show that they
can be viewed as refinements of the equivariant $L$-functions.

\begin{subsection}{The Stickelberger elements}\label{sec:theta}
\begin{definition}\label{def:thetadef}{\em{(\cite{g1},\cite{t1})}}
The Stickelberger element $\theta_{H'}=\theta_{L/K'}$ associated
to the extension $L/K'$ is the unique element of $\Z[H']$ such
that for each continuous character $\psi$ of $H'$,
\begin{equation}
\psi(\theta_{H'})=L_{S(K'),T(K')}(\psi,0).
\label{eq:thetadef}
\end{equation}
\end{definition}
For the existence of the Stickelberger element, see \cite{g1,t1}.

  In the case where $L/K$ is the constant $\Z_p$-extension, we
can relate $\theta_G$ to $\Theta_{\Gamma,S,T}$ in the following way.
First we note that the Galois group of the
constant $\Z_p$-extension over $k$ can be identified with some $H'$ such that
$G$ can be identified with
$\Gamma'\times H'$. This is actually the situation
discussed in Section \ref{sec:numerical}. We recall the notations
used there and in particular, we have $\varpi=|H'/H|$.
Let $\sigma'$ be the
Frobenius of $H'$. Then at every place $v\not\in S$, the Frobenius element
$[v]\in G$ can be expressed as the product
$\gamma_v'\cdot (\sigma')^{\deg(v)}$, $\gamma_v'\in\Gamma'$.
For every natural number $d$, there are only finitely many $v$ with
$\deg(v)=d$. From this we see that in the ring $\Z[\Gamma'][[\sigma']]$,
the infinite product
\begin{equation}\label{e:thetaprod}
\Pi=\prod_{v\in T} (1-N(v)\cdot[v])\prod_{v\not\in S}
(1-[v])^{-1}
\end{equation}
converges to a sum
$$1+\sum_{d=1}^{\infty} a_d'\cdot (\sigma')^d,$$
where each $a_d'$ is an element in the group ring $\Z[\Gamma']$.

\begin{lemma}\label{l:Thetatheta}
Suppose $L/K$ is the constant $\Z_p$-extension and $\sigma$ is the Frobenius
in $H=\Gal(L/K)$. Let notations be as the above. Then the followings hold.
\begin{enumerate}
\item We have $a_d'=0$, for almost all $d$, and also $\Pi=\theta_G$.
\item If $\yen(\theta_G)=\sum_d a_d\cdot (\sigma-1)^d,\; a_d\in\Z[\Gamma]$,
then $\Theta_{\Gamma,S,T}=\sum_d a_d\cdot (q^{-s}-1)^d$. Furthermore,
$a_d=0$, for $d=0,...,m$, if and only if $\theta_G\in I_H^m$.
In this case, we have
\begin{equation}\label{e:Thetatheta}
\Theta_{\Gamma,S,T}(0)^{(m)}=a_m=\Val_{\sigma,m,\breve{G}/H}(\yen([\theta_G]_{(m,H)}))
=\varpi^m\cdot \Val_{\sigma,m,G/H}([\theta_G]_{(m,H)}).
\end{equation}
\end{enumerate}
\end{lemma}
\begin{proof} We first apply $\yen$ to $\Pi$. Then we compare equations
(\ref{eq:ldef}), (\ref{e:Thetaprod}) and (\ref{e:thetaprod}) to see
that Part (1) and the first statement of Part (2) hold. The
rest is a consequence of Lemma \ref{l:yen}.
\end{proof}
\end{subsection}

\begin{subsection}{Twisted Stickelberger elements}\label{l:twist}

Suppose that $\{ H_{\gamma} \mid \gamma\in \Gamma \}$
are the $H$-cosets of $G$. We view $\Z[G]$ as the ring of integer valued
measures on $G$, and, for each $\gamma\in\Gamma$, let
$\Z[H_{\gamma}]$ be the set of integer valued measures on the
open subset $H_{\gamma}$. For a measure on $G$, we can restrict it to the
open subset $H_{\gamma}$, and this defines the restriction map
$
\res_{H_{\gamma}}:\Z[G]\longrightarrow \Z[H_{\gamma}].
$
Since $H_{\gamma}$ is also a closed subset of $G$, we can extend each
measure in $\Z[H_{\gamma}]$ to a
unique measure on $G$ vanishing outside $H_{\gamma}$. This
extending of measures induces the injective map
$\ext_{H_{\gamma}}:\Z[H_{\gamma}]\longrightarrow \Z[G]$.

\begin{definition}\label{def:retheta}
Define, for each $\gamma\in\Gamma$ and each $\xi\in\Z[G]$, the
$\gamma$-part of $\xi$ as $\xi_{\gamma}=\ext_{H_{\gamma}} \circ
\res_{H_{\gamma}} (\xi)$. We have $\xi=\sum_{\gamma} \xi_{\gamma}$.
\end{definition}

\begin{definition}\label{def:partheta}
For $\chi\in {\hat {\Gamma}}$, the $\chi$-twist homomorphism is
the ring homomorphism $[\chi]:\Z[\Gamma]\longrightarrow \O[\Gamma]$ which
sends $\sum a_{\gamma}\gamma$ to $ \sum
a_{\gamma}\chi(\gamma)\gamma$. Also, the
$\chi$-twist homomorphism from $\Z[G]$ to $\O[G]$ is the map
sending $\xi=\sum_{\gamma\in\Gamma} \xi_{\gamma}$ to
$\xi_{\chi}=\sum_{\gamma\in \Gamma}\chi(\gamma)\cdot
\xi_{\gamma}$, and we also let $[\chi]$ denote this
homomorphism.

\end{definition}
Thus we have the twisted Stickelberger elements $\theta_{\chi}=[\chi](\theta_{G})$, $\chi\in
\hat{\Gamma}$. We have the following commutative diagram.
\begin{equation}\label{e:chitwist}
\begin{array}{cccc}
[\chi]: & \Z[G] & \longrightarrow & \O[G]\\
{} & \downarrow & {} & \downarrow\\

[ \chi ]: & \Z[\Gamma] & \longrightarrow & \O[\Gamma],\\

\end{array}
\end{equation}
where two down-arrows are induced from the natural quotient map
$G\longrightarrow \Gamma$.

\end{subsection}

\end{section}

\begin{section}{The refined regulator maps}\label{sec:facts}
In this chapter we use the theory developed in Chapter \ref{sec:asaq}
to define the refined regulator maps.

\begin{subsection}{The global $\lambda$ map}
Let $\A_{K'}^*$ be the ideles group of $K'$ and
$\A_{K'}^*\longrightarrow H'$ be the norm residue map.
Recall the map $\delta_{H'}$ in (\ref{eq:del}). In a way similar to the one for constructing
the map ${\bar {\lambda}}_{v,\Gamma}$ in
(\ref{e:lvgamma}), for each place
$w$ of $K'$, composite the natural embeddings $U(K')\longrightarrow {K'}_w^*$ and
${K'}_w^*\longrightarrow \A_{K'}^*$ with $\delta_{H'}$ to form

\begin{equation}
\lambda_{w,H'}: U(K') \longrightarrow {K'}_w^*\longrightarrow
\A_{K'}^*\longrightarrow H'\stackrel{\delta_{H'}}{\longrightarrow}
I(H')/I(H')^2. \label{eq:unitmap}
\end{equation}
Following \cite{g2}, we define the global $\lambda$ map.
\begin{definition}\label{d:lambda}
Let
$$\lambda_{H'}:U(K')\longrightarrow X(K')\otimes I(H')/I(H')^2\subset
Y(K')\otimes I(H')/I(H')^2$$ be the homomorphism defined by
$$\lambda_{H'}(u)=\sum_{w\in S(K')} w\otimes\lambda_w(u),\;\mbox{for all }\;\;
u\in U(K').$$
\end{definition}
Note that we have, for $\gamma\in\Gamma'$,
\begin{equation}\label{e:equivlambda}
\lambda_{\gamma w,H}(\gamma u)=\gamma\lambda_{w,H}(u),
\end{equation}
and therefore
$\lambda_{H'}$ is a $\Gamma'$-equivariant homomorphism. Suppose that
$H'\subset H''=\Gal (L/K'')$ are two subgroups of $G$ and $w'$,
$w''$ are respectively places of $K'$ and $K''$ such that
$w''=w'\mid_{K''}$. We identify $Y(K'')$ with a sub-module of
$Y(K')$ by identifying $w''$ with the trace $\sum_{\sigma\in
H''/H'} \sigma w' \in Y(K')$. Also, we identify $I(H')/I(H')^2$
with a subgroup of $I(H'')/I(H'')^2$, using the diagram
$$
\begin{array}{ccc}
H' &\hookrightarrow& H''.\\
\downarrow &\circlearrowleft& \downarrow\\

I(H')/I(H')^2 &\hookrightarrow& I(H'')/I(H'')^2.\\

\end{array}
$$

 Then
both $Y(K')\otimes I(H')/I(H')^2$ and $Y(K'')\otimes
I(H'')/I(H'')^2$ can be viewed as sub-modules of $Y(K')\otimes
I(H'')/I(H'')^2$ in which we have, for $u\in U(K'')$,
\begin{equation}
\lambda_{H''}(u)=\lambda_{H'}(u).
\label{eq:compalambda}
\end{equation}
\end{subsection}

\begin{subsection}{The refined regulator maps}
\label{sec:regmap}
Let $\Lambda^n U\longrightarrow \Lambda^n (X(K)\otimes
I(H)/I(H)^2)$ be the homomorphism induced from
$\lambda_H$ and let
\newline
$\Lambda^n (X(K)\otimes
I(H)/I(H)^2)\longrightarrow \Lambda^n X(K)\otimes
I(H)^n/I(H)^{n+1}$ be the
$\Gamma$-equivariant homomorphism sending $(x_1\otimes
[h_1-1]_{(1)} )\wedge\dots\wedge (x_n\otimes [h_n-1]_{(1)})$ to
\newline
$(x_1\wedge\dots\wedge x_n)\otimes [(h_1-1)\cdot\dots\cdot
(h_n-1)]_{(n)}$. Then we define $\lambda_H^{(n)}$ as the composite
$$\lambda_H^{(n)}:\Lambda^n U\longrightarrow \Lambda^n (X(K)\otimes
I(H)/I(H)^2)\longrightarrow \Lambda^n X(K)\otimes
I(H)^n/I(H)^{n+1}.$$

Recall the notations in Chapter \ref{s:intro}. In particular, the map
$i^{(n)}$
%:\Lambda^n X(K)\longrightarrow \Lambda^n Y(K)$$
induces a map
$$i^{(n)}\otimes \id_H:\Lambda^n X(K)\otimes I(H)^n/I(H)^{n+1}
\longrightarrow \Lambda^n Y(K)\otimes I(H)^n/I(H)^{n+1}.$$
Also, for $\Psi\in \Lambda^n Y(K)^*$ the map $\iota(\Psi)$ induces a map
$$\iota(\Psi)\otimes \id_H: \Lambda^n Y(K)\otimes I(H)^n/I(H)^{n+1}\longrightarrow
\Z[\Gamma]\otimes I(H)^n/I(H)^{n+1}.$$
We then make the composition
$(\iota(\Psi)\otimes\id_H) \circ (i^{(n)}\otimes\id_H)\circ \lambda_H^{(n)}$ and extend it
linearly to form the map
$\mathcal{R}_{\Psi,H}^{\triangleright}: \Q\Lambda^n U\longrightarrow
\Q[\Gamma]\otimes I(H)^n/I(H)^{n+1}$.
\begin{definition}\label{d:regulatormap}
For an $\epsilon\in\Q\Lambda^n U$ such that $\mathcal{R}_{\Psi,\H}^{\triangleright}(\epsilon)$
is in $\Z[\Gamma]\otimes I(\H)^n/I(\H)^{n+1}$ and for the admissible Galois group $H$ with the natural
quotient map $Q_H:\H\longrightarrow H$
we define the refined regulator $\mathcal{R}_{\Psi,H}(\epsilon)$ as the image of
$\mathcal{R}_{\Psi,\H}^{\triangleright}(\epsilon)$ under the map
$$\Z[\Gamma]\otimes I(\H)^n/I(\H)^{n+1}\stackrel{\id\otimes Q_H}{\longrightarrow}
\Z[\Gamma]\otimes I(H)^n/I(H)^{n+1}\stackrel{\pounds_n}{\longrightarrow}
I_H^n/I_H^{n+1}.$$
\end{definition}
Similar to the map $\iota$ defined in Chapter \ref{s:intro}, we have
the map
$$\iota_H:\Lambda^n \Hom_{\Gamma}(U,\Z[\Gamma]\otimes
I(H)/I(H)^2)\longrightarrow \Hom_{\Gamma}(\Lambda^n U, \Z[\Gamma]\otimes
I(H)^n/I(H)^{n+1}),$$
such that if $\phi_1,...\phi_n\in \Hom_{\Gamma}(U,\Z[\Gamma]\otimes
I(H)/I(H)^2)$ and $u_1,...,u_n\in U$, then
\begin{equation}\label{e:iotaH}
\iota_H(\phi_1\wedge\dots \wedge \phi_n)(u_1\wedge\dots\wedge
u_n)=\det(\phi_i(u_j)),
\end{equation}
where the determinant is computed by using the multiplication
\begin{equation}\label{e:mutipletensor}
\begin{array}{rcl}
\Z[\Gamma]\otimes I(H)^i/I(H)^{i+1}\times
\Z[\Gamma]\otimes I(H)^j/I(H)^{j+1}
&\longrightarrow
&\Z[\Gamma]\otimes I(H)^{i+j}/I(H)^{i+j+1}\\
(a\otimes [\xi_i]_{(i)}, b\otimes [\xi_j]_{(j)})
& \mapsto &
ab\otimes [\xi_i\xi_j]_{(i+j)}.\\
\end{array}
\end{equation}
If $\psi_i\in Y(K)^*$, then $\mathcal{R}_{\psi_i,H}^{\triangleright }\in
\Hom_{\Gamma}(U,\Z[\Gamma]\otimes
I(H)/I(H)^2)$, and we have
\begin{equation}\label{e:itotaetaH}
\mathcal{R}_{\psi_1\wedge\dots\wedge\psi_n,H}^{\triangleright}=\iota_H
(\mathcal{R}_{\psi_1,H}^{\triangleright}
\wedge\dots\wedge \mathcal{R}_{\psi_n,H}^{\triangleright}).
\end{equation}

Next, we will study the relation between $\mathcal{R}_{\Psi}$ and $\mathcal{R}_{\Psi,H}$.
Our first goal is to show that $\mathcal{R}_{\Psi,H}(\epsilon)$ is defined for every
$\epsilon\in\Z_{(p)}\Lambda^n_0 U$.

Suppose $H\simeq \Z_p^d$ and $\{t_1,...,t_d\}$ is a basis. Then $\Z[\Gamma]\otimes I(H)/I(H)^2$
can be identified with the direct sum $\Z_p[\Gamma]t_1+\dots +\Z_p[\Gamma]t_d$.
Let $\psi_1,...,\psi_n\in Y(K)^*$. From the above construction, each
$\mathcal{R}_{\psi_i,H}^{\triangleright}$ is a $\Gamma$-equivariant map from
$\Q U$ to $\Q_p[\Gamma]t_1+\dots +\Q_p[\Gamma]t_n$ and we
have $\mathcal{R}_{\psi_i,H}^{\triangleright}=\oplus_{j=1}^d \psi_{ij}t_j$ where each $\psi_{ij}$ is
an element of $\Z_{(p)} U^*$. From the isomorphism (\ref{eq:mapd}) and Lemma \ref{le:basechange},
we see that $I(H)^n/I(H)^{n+1}$ is the $n$th symmetric tensor of $I(H)/I(H)^2$ and
$\Z[\Gamma]\otimes I(H)^n/I(H)^{n+1}$ can be
identified with the direct sum $\bigoplus_{n_1+\dots +n_d=n}
\Z_p[\Gamma] t_1^{n_1}\cdot\dots\cdot t_d^{n_d}$.
If $\Psi=\psi_1\wedge\dots\wedge \psi_n$, then we have
$\mathcal{R}_{\Psi,H}^{\triangleright}=
\oplus_{n_1+\dots n_d=n} \Psi_{n_1,...,n_d}t_1^{n_1}\cdot\dots\cdot t_d^{n_d}$ where
each $\Psi_{n_1,...,n_d}$ is a $\Gamma$-equivariant map from $\Q\Lambda^n U$ to $\Q[\Gamma]$.
In fact, if $\Xi_{n_1,...,n_d}$ is the set consisting
of all maps
\newline
$\xi:\{1,...,n\}\longrightarrow \{1,...,d\}$ such that
$|\xi^{-1}(i)|=n_i$,
then from (\ref{e:itotaetaH}) we see that
$$\Psi_{n_1,...,n_d}=\sum_{\xi\in\Xi_{n_1,...n_d}} \iota(\psi_{1 \xi(1)}\wedge\dots\wedge
\psi_{n \xi(n)}).$$
For each $\epsilon\in\Z_{(p)}\Lambda^n_0 U$, we have $ \iota(\psi_{1 \xi(1)}\wedge\dots\wedge
\psi_{n \xi(n)})(\epsilon)\in\Z_p[\Gamma]$, hence we must have
$\mathcal{R}_{\Psi,H}^{\triangleright}(\epsilon)\in\Z[\Gamma]\otimes I(H)^n/I(H)^{n+1}$.
Furthermore, since $\H$ is the projective limit of those $H$ which are finite free over $\Z_p$,
we also have $\mathcal{R}_{\Psi,\H}^{\triangleright}(\epsilon)\in\Z[\Gamma]\otimes
I(\H)^n/I(\H)^{n+1}$ for $\epsilon\in\Z_{(p)}\Lambda^n_0 U$. The following lemma is proved.
\begin{lemma}\label{l:lambda0}
If $\Psi\in \Lambda^n Y(K)^*$, then the refined regulator $\mathcal{R}_{\Psi,H}(\epsilon)$ is defined for every
$\epsilon\in\Z_{(p)}\Lambda^n_0 U$.
\end{lemma}
Now we compare the refined regulator map with the old regulator map.
\begin{lemma}\label{l:regrelation}
Suppose $L/K$ is the constant $\Z_p$-extension and $\sigma$ is the Frobenius in $H=\Gal(L/K)$.
Then for every $\Psi\in \Lambda^n Y(K)^*$ and every $\epsilon\in \Lambda_0^n U$
\begin{equation}\label{e:regrelation}
\mathcal{R}_{\Psi}(\epsilon)= \Val_{\sigma,n}(\yen(\mathcal{R}_{\Psi,H}
(\epsilon))).
\end{equation}
\end{lemma}
\begin{proof}
It is easy to see that $\Val_{\sigma,1}(\yen(\pounds_1(\lambda_{w,H})))=\lambda_w$,
and from this we see that
$\Val_{\sigma,1}(\yen(\pounds_1(\mathcal{R}_{\psi_i,H}^{\triangleright}(u_j))))
=\mathcal{R}_{\psi_i}(u_j)$. This proves the lemma for the $n=1$ case. The general case
is proved by using (\ref{e:etaiota}) and (\ref{e:itotaetaH}).
\end{proof}
\end{subsection}

\begin{subsection}{The twisted regulators}
Recall the $\chi$-twist homomorphisms defined in Definition
\ref{def:partheta}.

\begin{lemma}\label{le:twisthom} We have $$[\chi](I_H^n)\subset
I_{{\O},H}^n,$$
\end{lemma}
\begin{proof}
This is due to (\ref{e:chitwist}).
\end{proof}
By the abuse of notations, we also use $[\chi]$ to denote the induced homomorphism
$I_H^n/I_H^{n+1}\longrightarrow
I_{{\O},H}^n/I_{{\O},H}^{n+1}.$

\begin{definition}\label{def:twistregulator}
For $\chi\in {\hat {\Gamma}}$, $\Psi\in \Lambda^n Y(K)^*$ and
$\epsilon\in\Z_{(p)}\Lambda_0^n U$
define the $\chi$-twisted regulator
$$\mathcal{R}_{\Psi, \chi,H}(\epsilon)=[\chi](\mathcal{R}_{\Psi,H}(\epsilon))
\in I_{\O,H}^n/I_{\O,H}^{n+1}.$$

\end{definition}

\end{subsection}
\begin{subsection}{The universal property}
\label{sec:unique} In this section we study some special properties of the refined
regulator map. We will show the injectivity as well as a universal property which says that
every set of $n$ homomorphisms $\phi_1,...,\phi_n\in\Hom_{\Gamma}(U,\Z_p[\Gamma])$ can be obtained
from the refined regulator map.

%The main result of this section is Lemma \ref{le:uniqueintegral}
%which will imply the uniqueness of the element
%$\epsilon\in\Lambda_{S,T}$ in Theorem \ref{thm:stark1} and Theorem
%\ref{thm:stark2}. It will
%also show that if this $\epsilon$ exists then it will be automatically
%in $\Lambda_{\H}^n U$. To prove this lemma, we need to use the following
%two lemma.  To simplify the notations, through the map
%$\pounds_n$, we will identify
%$\Z[\Gamma]\otimes I(H)^n/I(H)^{n+1}$ with $I_{H}^n/I_{H}^{n+1}$. Then
%the following is obvious.

Let $w_1,..,w_n,\eta$ and $w_i^*,...,w_n^*$ be those defined in Chapter \ref{s:intro}.
Let $S'$ be a finite set
of places of $k$ such that $S\cap S'=\emptyset$. For each $v\in S'$ we arbitrarily
choose a place $w$ of $K$ sitting over $v$. Let $S_K'$ be the collection of these
chosen places.
\begin{lemma}\label{le:inject} For a place $w'\in S_K'$ let $[w']\in \H$ be the Frobenius element
at $w'$.
If there are $u_1,...,u_n\in U$ and $a_{w'}\in\Z_p$ for each $w'\in S_K'$
such that in $I(\H)/I(\H)^2$ the sum $\sum_{i=1}^n
\lambda_{w_i,\H}(u_i)+\sum_{w'\in S_K'}a_{w'}\delta_{\H}([w'])
$ is divisible by $p$, then in
$U$ every $u_i$ is divisible by $p$.
\end{lemma}
\begin{proof} We consider the maximal pro-$p$ abelian extension ${\tilde L}/K$
unramified outside $S(K)$. Denote ${\tilde{\G}}=\Gal({\tilde L}/k)$ and
${\tilde {\H}}=\Gal({\tilde L}/K)$. Then $\Gamma={\tilde{\G}}/{\tilde {\H}}$
and $[{\tilde{\G}},{\tilde{\G}}]\subset {\tilde {\H}}$. In particular, $\Gamma$
acts on ${\tilde{\H}}$ through conjugate. Let
$$\mathbb{K}=[{\tilde{\G}},{\tilde{\H}}]=
\sum_{\gamma\in\Gamma}(1-\gamma){\tilde{\H}},$$
and denote $\G'={\tilde{\G}}/{\mathbb{K}}$,
$\H'={\tilde{\H}}/{\mathbb{K}}$. Then, $\Gamma=\G'/\H'$ and $\H=\H'/[\G',\G']$. As $\Gamma$ acts
trivially on $\H'$, we have
\begin{equation}\label{eq:commuta}
[\G',\H']=\{\id\}.
\end{equation}
Let $N=|\Gamma |$. Then for $f,g\in\G'$, $fg^Nf^{-1}=g^N$ and hence
$fgf^{-1}=gh$ for some $h\in\H'$ such that $h^N=id$. By (\ref{eq:commuta}),
the abelian group $[\G',\G']$ is generated by $N^2$ elements whose orders are divisors of $N$.
Therefore, $[\G',\G']$ is finite and so is the index
$C:=[[{\tilde{\G}},{\tilde{\G}}]:\mathbb{K}]$.
Let $p^m$ be the maximal $p$-power divisor of $C$.
\par
Let $z_i=(z_{i,w})_w \in \A_K^*$ be the idele such that
$z_{i,w}=u_i$ if $w=w_i$ and $z_{i,w}=1$ if $w\not=w_i$. If
$\sum_{i=1}^n \lambda_{w_i,\H}(u_i)+\sum_{w'\in S_K'} a_{w'}\delta_{\H}([w'])$
is divisible by $p$ in
$I(\H)/I(\H)^2$, then by Class Field Theory, we can find
$\alpha\in K^*$, $\beta\in \prod_{w\notin S(K)} \O_w^*$, $a\in \A_K^*$ and
$b_{\gamma}\in \A_K^*$ for each $\gamma\in \Gamma$ such that
$$\prod_{i=1}^n z_{i}^C=a^{p^{m+1}}\cdot \alpha\cdot \beta
\cdot\prod_{w'\in S_K'}\pi_{w'}^{Ca_{w'}}
\cdot\prod_{\gamma\in\Gamma}\gamma b_{\gamma} /b_{\gamma}.$$

Note
that since $\{ v_1,...,v_n\}=S_0\varsubsetneq S$, there is a place $v_0\in S
\backslash S_0$. Taking the norm $N_{K/k}$, we see that at $v_0$ the norm
$N_{K/k}(\alpha)$ is locally a $p^{m+1}$th power. By Lemma
\ref{le:localleo}, the element $N_{K/k}(\alpha)$ is a $p^{m+1}$th power
in $k^*$. This implies that each $N_{K/k}(z_{i})$ is a $p$th power
idele. Since $z_{i}$ is trivial away from $w_i$ and $v_i$ splits
completely in $K$, $u_{i}=N_{K/k}(z_{i})$ itself is a $p$th power
in $K_{w_i}^*$. Again, Lemma \ref{le:localleo} implies that $u_i$
is a $p$th power in $K^*$ and hence a $p$th power in $U$.
\end{proof}
Let $U_i=\Z_p\otimes U$, for $i=1,...,n$. We extend each
$\lambda_{w_i,\H}$ linearly to a map
from $U_i$ to $I(\H)/I(\H)^2$ and form the sum
$$\lambda_{S_0,\H}:=\sum_{i=1}^n \lambda_{w_i,\H}: U_1\oplus \dots\oplus U_n
\longrightarrow  I(\H)/I(\H)^2.$$
Let $\mathbb{U}$ be the image of $\lambda_{S_0,\H}$ and
denote by $\mathbb{W}$ the $\Z_p$-sub-module of $ I(H)/I(H)^2$
generated by the set $\{\delta_{\H}([w']) \mid \; w'\in S_K'\}$.
In view of Lemma \ref{le:inject}, we have the following.

\begin{lemma}\label{le:inject2}
The map $\lambda_{S_0,\H}$ is injective and its image, denoted as $\mathbb{U}$,
is a direct summand of $I(\H)/I(\H)^2$ as a $\Z_p$-sub-modules. Furthermore,
we have $\mathbb{U}  \cap \mathbb{W} =\{0\}$ and $\mathbb{U}+\mathbb{W}/\mathbb{W}$
is a direct summand of $(I(\H)/I(\H)^2)/\mathbb{W}$.
\end{lemma}
%\begin{proof}
%It is easy to see that $\lambda_{S_0,\H}$ is injective
%and $\mathbb{U}  \cap \mathbb{W} =\{0\}$.
%To simplify the notation, we use the isomorphism $\delta_{\H}$ to identify $I(\H)/I(H)^2$ with
%$\H$. Then $\mathbb{U}$ and $\mathbb{W}$
%are identified with sub-modules of $\H$.

%Let $H$ be an admissible Galois group finite
%free over $\Z_p$. Denote by $U_H$ and $W_H$ the image of $\mathbb{U}$ and $\mathbb{W}$
%under the natural projection $\H\longrightarrow H$. we say that the statement $(P)$ holds if
%for any $z\in H$ with $pz\in U_H+W_H$
%there is a $u\in U_H$ to make $p(z-u)\in W_H$.
%If the statement $(P)$ holds for $H$, then $U_H$ is a direct summand of $H$ and $U_H+W_H/W_H$
%is a direct summand of $H/W_H$.
%Suppose $L_1/k$ and $L_2/k$ admissible extensions with $L_1\subset L_2$
%such that the Galois groups $H_1=\Gal(L_1/K)$ and $H_2=\Gal(L_2/K)$ both finite free over
%$\Z_p$. It is easy to see that if the statement $(P)$ holds for $H_1$
%then it also holds for $H_2$. Thus it is enough to show that the statement $(P)$ holds
%for some $H$.

%\end{proof}
As in Chapter \ref{s:intro}, if $M$ is a $\Z[\Gamma]$-module then there is a one-one
correspondence $\phi \leftrightarrow\phi^{(id)}$  between
$\Hom_{\Gamma}(M,I_{H}/I_{H}^2)$ and
$\Hom_{\Z}(M,I(H)/I(H)^2)$
such that
$$\phi(m)=\sum_{\gamma\in\Gamma} \gamma\otimes \phi^{(id)} ( \gamma^{-1} (m)),$$
for every $m\in M$. A similar correspondence holds if $M$ is a $\Z_p[\Gamma]$-module.

Recall that $\mathcal{R}_{w_i^*,\H}^{\vartriangleright}$ is formed by the
linear extending of the composition
$(w_i^*\otimes\id_{\H})\circ (i\otimes \id_{\H})\circ \lambda_{\H}$. Using the above
notations, we easy find that
\begin{equation}\label{e:(id)}
((w_i^*\otimes\id_{\H})\circ (i\otimes \id_{\H})\circ \lambda_{\H})^{(\id)}=\lambda_{w_i,\H}.
\end{equation}

\begin{corollary}\label{c:(id)}
The maps $\mathcal{R}_{w_1^*,\H}^{\vartriangleright},...,
\mathcal{R}_{w_n^*,\H}^{\vartriangleright}\in \Hom_{\Gamma}(\Q U,\Q[\Gamma]\otimes I(\H)/I(\H)^2)$
are linearly independent over $\Q_p$.
\end{corollary}
\begin{proof}
Lemma \ref{le:inject2} says $(\mathcal{R}_{w_1^*,\H}^{\vartriangleright})^{(id)},...,
(\mathcal{R}_{w_n^*,\H}^{\vartriangleright})^{(id)}$ are linearly independent.
\end{proof}
\begin{corollary}\label{c:projection}
Let notations be as in {\em{Lemma} \ref{le:inject2}}.
For every $\phi_1,...,\phi_n\in \Hom_{\Gamma}(U,\Z[\Gamma]\otimes I(H)/I(H)^2)$ there is
a $\Z_p$-morphism $e:\H\longrightarrow H$ such that $\phi_i=
\delta_H \circ e\circ \delta_{\H}^{-1}\circ \mathcal{R}_{w_i^*,\H}^{\vartriangleright}$
for every $i$.
Furthermore, $e$ can be chosen such that
$e(\mathbb{W})=0$.
\end{corollary}
\begin{proof}
Lemma \ref{le:inject2} says that $e$ can be chosen to satisfy
$(\phi_i)^{(id)}=
\delta_H \circ e\circ \delta_{\H}^{-1}\circ (\mathcal{R}_{w_i^*,\H}^{\vartriangleright})^{(id)}$
for every $i$ and $e(\mathbb{W})=0$.
\end{proof}

We make the identification $\Z_p=I(H_0)/I(H_0)^2=H_0:=
\Gal(K\F_{q^{p^{\infty}}}/K)$ and view $\Z$ as subgroup of
$I(H_0)/I(H_0)^2$. Then
$U^*=\Hom_{\Gamma}(U,\Z[\Gamma])$ is identified as a subgroup of
$\Hom_{\Gamma}(U,\Z[\Gamma]\otimes I(H_0)/I(H_0)^2)$.

Over $k$ choose a finite number of places not in $S$ such that with
respect to
the abelian extension $K/k$ the decomposition subgroups at these
places generate the Galois group $\Gamma$. Let $S'$ denote the set form
by these places.
Suppose $\Phi=\iota(\phi_1\wedge\dots\wedge\phi_n)\not=0$ with
$\phi_1,...,\phi_n\in U^*\subset \Hom_{\Gamma}(U,\Z[\Gamma]\otimes I(H_0)/I(H_0)^2)$.
By Corollary \ref{c:projection},
there is a $\Z_p$-morphism $e:\H\longrightarrow H_0$ such that
\begin{equation}\label{e:phiri}
\phi_i=
\delta_{H_0} \circ e\circ \delta_{\H}^{-1}\circ
\mathcal{R}_{w_i^*,\H}^{\vartriangleright},\; \text{for}\; i=1,...,n.
\end{equation}
Since $\Phi\not=0$, the co-kernel of the morphism $e$ must be finite.
Let $c$ be its order, and
let $L_1$ be the fixed field of the kernel of $e$ and denote
$H_1=\Gal(L_1/K)$, $G_1=\Gal(L_1/k)$. We have $e=j\circ Q$
where $Q:\H\longrightarrow H_1$ is the natural projection and
$j:H_1\stackrel{\sim}{\longrightarrow} cH_0 \hookrightarrow H_0$.
Furthermore, Corollary \ref{c:projection} says that $e$ can be chosen such
that every place of $K$ sitting over $S'$ splits completely over $L_1/K$.
Thus, over $L_1/k$ the decomposition subgroup at each place in $S'$
is a finite subgroup of $G_1$ and these decomposition groups generate
a finite group which, under the natural projection $G_1\longrightarrow
G_1/H_1=\Gamma$, is isomorphic to $\Gamma$. This means that
$G_1$ is the direct product $\Gamma\times H_1$.

Let $\sigma\in H_1$
be the generator such that $j(\sigma)$ is $c$ times the Frobenius in $H_0$,
and let $t=\sigma-1\in\Z[H_1]$. Then we have
\begin{equation}\label{e:formal}
F_p[G_1]=F_p[\Gamma][[t]]
\end{equation}
and by
Lemma \ref{le:torfree}, the $n$th relative augmentation ideal
$I_{F_p,H_1}^n$ is just the principal ideal $(t^n)$. If $\epsilon_1$
is an element in $\Q\Lambda^n U$, then
$\pounds_n(\mathcal{R}_{\eta,H_1}^{\vartriangleright}(\epsilon_1))
=a_nt^n$ for some
$a_n\in F_p[\Gamma]$ and we have
$$
\Phi(\epsilon_1)=\iota(\phi_1\wedge\dots\wedge \phi_n)(\epsilon_1)=
a_n=\Val_{\sigma,n,G/H_1}(\pounds_n(\mathcal{R}_{\eta,H_1}^{\vartriangleright}
(\epsilon_1))).
$$
In particular, if $\mathcal{R}_{\eta,\H}(\epsilon_1)$
is defined, then $\Phi(\epsilon_1)\in \Z_p[\Gamma]$,
and hence $\epsilon_1$ is
an element of $\Z_{(p)}\Lambda_0^n U$. Thus we have prove the following
corollary.
\begin{corollary}\label{c:lambda0H}
If $\mathcal{R}_{\eta,\H}(\epsilon)$ is defined, then $\epsilon$ is indeed
an element of $\Z_{(p)}\Lambda_0^n U$. If notations are as above,
then we have
\begin{equation}\label{e:phirelation}
\Phi(\epsilon)=\Val_{\sigma,n,G/H_1}(\mathcal{R}_{\eta,H_1}
(\epsilon)).
\end{equation}
\end{corollary}
%\begin{proof}
%Suppose that $\mathcal{R}_{\eta,\H}(\epsilon)$ is defined. We need to show
%that for $\phi_1,...,\phi_n\in \Hom_{\Gamma}(U,\Z[\Gamma])$,
%the image of $\epsilon$ under the map $\iota(\phi_1\wedge\dots
%\wedge \phi_n)$ is in $\Z_{(p)}[\Gamma]$. We prove this by taking
%$\Z\subset H=:\Gal(K\F_{q^{p^{\infty}}}/K)$ and view each $\phi_i$ as
%in $ \Hom_{\Gamma}(U,\Z[\Gamma]\otimes I(H)/I(H)^2)$. Then we use
%(\ref{e:itotaetaH}) and Corollary \ref{c:projection} to complete the proof.

%\end{proof}
\begin{lemma}\label{le:uniqueintegral}
Suppose that $n$ and $\eta$ are as those in {\em {Conjecture
\ref{c:r}}}.
%, and $ \lambda_{v_1},...,\lambda_{v_n} \in
%\Hom_{\Gamma}(U,I_{\H}/I_{\H}^2)$ are as those in {\em {Definition
%\ref{d:wstar'}}}.
Then the followings are true.
\begin{enumerate}
\item The map $\mathcal{R}_{\eta,\H}^{\triangleright}:\Q \Lambda^n U \longrightarrow
\Q  I_{\H}^n/I_{\H}^{n+1}$ is injective.
\item For each $\chi\in {\hat {\Gamma}}$, the map
$\mathcal{R}_{\eta,\chi,\H}^{\triangleright}=:[\chi]\circ
\mathcal{R}_{\eta,\H}^{\triangleright}:\Q \Lambda^n U \longrightarrow
F I_{\H}^n/I_{\H}^{n+1}$ is injective.
%\item For $\chi\in {\hat {\Gamma}}$, the map
%$$\Q\cdot\Lambda^{n(\chi)} U(K_{\chi}) \stackrel{\mathcal{R}_{\eta,\chi,\H}}{\longrightarrow}
%F\cdot I_{\H}^{n(\chi)}/I_{\H}^{n(\chi)+1}\stackrel{t}{\longrightarrow}
%F\cdot I(\H)^{n(\chi)}/I_(\H)^{n(\chi)+1}$$
%is injective.
%\item
%\item We have
%$$\Lambda_{\H}^n U=\Lambda_{\mathcal{R}_{\eta,\H}}^n U.$$
\end{enumerate}
\end{lemma}

\begin{proof}
The $\chi$-twist map $[\chi]:\Q  I_{\H}^n/I_{\H}^{n+1}\longrightarrow
F I_{\H}^n/I_{\H}^{n+1}$ is injective, and hence Part (2) is a
consequence of Part (1).

To prove Part (1), We first apply Corollary \ref{c:(id)} and
choose an admissible $H_0$
such that $H_0\simeq\Z_p^d$ for some $d$ and
the subspaces $\mathcal{R}_{w_1^*,H_0}^{\vartriangleright}(U),...,
\mathcal{R}_{w_n^*,H_0}^{\vartriangleright}(U)
\subset \Q[\Gamma]\otimes I(H_0)/I(H_0)^2$
are linearly independent over $\Q_p$.
%Here, as usual, $Q_*$ is induced
%from the natural quotient map $Q:\H\longrightarrow H_0$.
To simplify the notations, put $W={\bar{\Q}_p}\otimes U$,
$V={\bar{\Q}_p}\otimes_{\Z_p}I_{H_0}/I_{H_0}^2$ and
linearly extend each $\mathcal{R}_{w_i^*,H_0}^{\vartriangleright}$
to the map $\mathcal{R}_i:W\longrightarrow V$.

From Equation (\ref{e:itotaetaH}), we see that it is
enough to show that the map
$$\iota_{H_0}(\mathcal{R}_1\wedge\dots\wedge\mathcal{R}_n):\Lambda^n W
\longrightarrow {\bar{\Q}_p}\otimes_{\Z_p} I_{H_0}^n/I_{H_0}^{n+1}$$
is injective.
%have
%$$\mathcal{R}_{\eta,\H}^{\triangleright}=\iota_{\H}(
%\mathcal{R}_{w_1^*,\H}^{\triangleright}
%\wedge\dots\wedge \mathcal{R}_{w_n^*,\H}^{\triangleright}).$$

 By Lemma \ref{le:torfree} and Equation (\ref{eq:ideal}),
in the category of ${\bar{\Q}_p}[\Gamma]$-module,  ${\bar{\Q}_p}\otimes_{\Z_p}
I_{H_0}^n/I_{H_0}^{n+1}$ is nothing but the $n$th symmetric tensor
of $V$. Without loss of generality, we can assume that
$V=\mathcal{R}_1(W)\oplus\dots\oplus\mathcal{R}_n(W)$. If
$$W=W_1+\dots+W_m$$
is the decomposition of $W$ into irreducible
${\bar{\Q}_p}[\Gamma]$-modules, then $\Lambda^n W$ is decomposed into the
direct sum $\bigoplus_A W_A$, where, associated to each
$A=\{i_1,...,i_n\}\subset \{1,...,m\}$ such that $|A|=n$, $W_A$ is
the exterior tensor of $W_{i_1},..., W_{i_n}$. Also, ${\bar{\Q}_p}\otimes_{\Z_p}
I_{H_0}^n/I_{H_0}^{n+1}$ is decomposed into the direct sum
$\bigoplus_{A,\sigma} V_{A,\sigma}$, where associated to each pair
$(A,\sigma)$ with $A$ as above and $\sigma\in S_n$, the symmetric
group of $n$ elements, $V_{A,\sigma}$ is the symmetric tensor of
$\mathcal{R}_{1}(W_{\sigma(i_1)}),..., \mathcal{R}_{i_n}(W_{\sigma(i_n)})$.
By (\ref{e:iotaH}), the homomorphism
$\iota_{H_0}(\mathcal{R}_1\wedge\dots\wedge\mathcal{R}_n)$
is injective on each $W_A$ and it sends $W_A$ into $\sum_{\sigma\in
S_n} V_{A,\sigma}$. Therefore, it is
injective on $\Lambda^n W$.

\end{proof}

\end{subsection}
\end{section}

\begin{section}
{The main theorem}
In this chapter, we state our main theorem and show that it implies Theorem
\ref{th:rp} and Theorem \ref{th:burnsp}. We also state a twisted version of it.

\begin{subsection}{The main theorem}
\label{sec:mainthm}

\begin{theorem}\label{thm:stark1}
Let notations be as in {\em{Corollary \ref{c:r}}}. Then for every admissible $H$,
the Stickelberger element $\theta_G\in I_H^n$. Furthermore,
there is a unique
$\epsilon\in \Z_{(p)}\Lambda_{S,T}^n$ such that for every admissible $H$
\begin{equation}
{\mathcal{R}}_{\eta,H}(\epsilon)=[\theta_G]_{(n,H)}
\label{eq:stark}
\end{equation}
\end{theorem}
Note that the uniqueness of $\epsilon$ follows from (\ref{eq:stark}) and
Lemma \ref{le:uniqueintegral}.

Now we show that this main theorem implies Theorem
\ref{th:rp}.
\begin{proof}(of Theorem \ref{th:rp})
Let $H=\Gal(K\F_{q^{p^{\infty}}}/K)$ and let $\sigma$ be the Frobenius,
and apply $\Val_{\sigma,n,G/H}\circ\yen$ to both side of (\ref{eq:stark}).
Then we use Lemma \ref{l:Thetatheta} and Lemma \ref{l:regrelation}.
\end{proof}
\end{subsection}

\begin{subsection}{The Conjecture of Gross}
\label{sec:gross}
To complete the proof of Theorem \ref{th:burnsp} we need to use some
results concerning a conjecture of Gross, which will
be discussed in this section.
This conjecture can be viewed as a refinement of the class number
formula in which ${\det}_{H'}$, an refined regulator of Gross, is involved.

In \cite{g1}, this refined regulator is defined
as an element in $I^n/I^{n+1}$.  We instead choose to
adopt  Tate's definition \cite{t2} and define the refined regulator  as an element in
the group ring. Here we describe Tate's definition of the refined regulator.

For a place $w\in S(K')$ let
$\lambda_{w,H'}:U(K')\longrightarrow I(H')/I(H')^2$ be the map in
(\ref{eq:unitmap}), and let $\delta_{H'}$ be the map in
(\ref{eq:del}). Suppose $w_1',...,w_{r_{K'}}'$ are distinct places
in $S(K')$ and $u_1,...,u_{r_{K'}}$ is a $\Z$-basis of $U(K')$, and assume
that the ordering of them are chosen such that the classical regulator
form by them is positive. Then the refined regulator of Gross is defined as
\begin{equation}\label{e:greg}
{\det}_{H'}=\det_{1\leq i, j \leq r_{K'}}
(\delta_{H'}^{-1}(\lambda_{w_i',H'}(u_j))-1)\in
I(H')^{r_{K'}}.
\end{equation}

For each $K'$, recall that $h_{K',S(K'),T(K')}$ is the modified class number
defined in (\ref{e:hkst}).

\begin{conjecture}\label{c:gross} \em{(Gross)}
We have $\theta_{H'}\in I(H')^{r_{K'}}$ and
$$\theta_{H'}\equiv h_{K',S(K'),T(K')}\cdot {\rm det}_{H'}
\pmod{I(H')^{r_{K'}+1}}.$$
\end{conjecture}
For evidences and related discussions of this conjecture, see,
for examples, \cite{a91,a03,b02,b05,bl04,d95,h,h04,l97,l02,l04,rd02,t2,t04,
tn,tn2,y89}.
We will need the following result from \cite{tn}.

\begin{theorem}\label{th:gross}
We have $\theta_{H'}\in I_p(H')^{r_{K'}}$ and
\begin{equation}\label{e:grossp}
\theta_{H'}\equiv h_{K',S(K'),T(K')}\cdot {\rm det}_{H'}
\pmod{I_p(H')^{r_{K'}+1}}.
\end{equation}
\end{theorem}
Note that in \cite{tn}, this theorem is proved only for the case where
$H'$ is a pro-$p$ group. But, since $\theta_{H'}$ is known to be
in $I(H')$ and if $H'=H_p'\times H_0'$ is the decomposition into
the direct product of the pro-$p$ part, $H_p'$, and the non-$p$ part, $H_0'$,
then for each $m\geq 1$,
$$I_p(H')^m/I_p(H')^{m+1}=
I_p(H_p')^m/I_p(H_p')^{m+1}\times I_p(H_0')^m/I_p(H_0')^{m+1},$$
with $I_p(H_0')^m/I_p(H_0')^{m+1}=0$. Therefore,
the theorem holds for general $H'$.

We are ready to show that Theorem \ref{thm:stark1} implies
Theorem \ref{th:burnsp}.

\begin{proof} (of Theorem \ref{th:burnsp})
We recall the notations in the proof of Corollary \ref{c:lambda0H}.
Thus the admissible Galois groups $H_0$ and $H_1$ are chosen and we have
$G_1=\Gamma\times H_1$, also $\sigma\in H_1$ is a
generator
and with $t=\sigma-1\in\Z[H_1]$ we have
$$
\Z_p[G_1]=\Z_p[\Gamma][[t]]
$$
such that the $n$th relative augmentation ideal
$I_{p,H_1}^n$ is just the principal ideal $(t^n)$.
The first part of Theorem
\ref{thm:stark1} say that $\theta_{G_1}\in I_{p,H_1}^n$ and hence
\begin{equation}\label{e:lfirstterm}
\theta_{G_1}=a_nt^n+\dots+a_{i}t^{i}+\dots, \; a_i\in \Z_p[\Gamma],
\end{equation}
from which, we get
$a_n=\Val_{\sigma,n,G/H_1}([\theta_{G_1}]_{(n,H_1)})$.
Now Corollary \ref{c:lambda0H} and the second part of
Theorem \ref{thm:stark1} says that

\begin{equation}\label{e:rfirstterm}
\Phi(\epsilon)=\Val_{\sigma,n,G/H_1}( \mathcal{R}_{\eta,H_1} (\epsilon) )
=\Val_{\sigma,n,G/H_1}( [ \theta_{G_1} ]_{ (n,H_1) } )=a_n.
\end{equation}

If $v_i\in S_0$, then it splits completely over $K$ and hence for $u\in U(k)$ the image
$\lambda_{v_i,G_1}(u)\in I(H_1)/I(H_1)^2\subset I(G_1)/I(G_1)^2$.
This implies that ${\det}_{G_1}$ is in
$I(G_1)^{r-n}I(H_1)^n$ and hence
\begin{equation}\label{e:rgrossterm}
{\det}_{G_1}=b_nt^n+\dots+b_it^i+\cdots.
\end{equation}

Since $G_1$ is the direct product of $\Gamma$ and $H_1$,
the augmentation ideal $I_p(G_1)$ is generated by
$t$ and $I_p(\Gamma)$. Therefore, an element
$$\xi=c_0+c_1t+\dots+c_it^i+\dots \in\Z_p[G_1]$$
is in $I_p(G_1)^m$ if and only if $c_i\in I_p(\Gamma)^{m-i}$ for every $i\leq m$.
This together with (\ref{e:lfirstterm}), (\ref{e:rgrossterm})
and Theorem \ref{th:gross} implies that for $i=n,...,r_k$

\begin{equation}\label{e:aandb}
a_i\equiv b_i\pmod{I_p(\Gamma)^{r_k-i+1}}.
\end{equation}

From (\ref{e:phiri}), we see that for $i=1,...,n$ and $u\in U(k)$
$$\delta_{G_1}^{-1}(\lambda_{v_i,G_1}(u))-1=\phi_i^{(id)}(u)t.$$
For $i=n+1,...,r_k$ and $u\in U(k)$,
we have
$$\delta_{G_1}^{-1}(\lambda_{v_i,G_1}(u))={\bar {\lambda}}_{v_i,\Gamma}(u)\cdot\sigma^{\alpha},\;
\text{for some}\;\alpha\in\Z_p,$$
and hence
$$\delta_{G_1}^{-1}(\lambda_{v_i,G_1}(u))-1\equiv
{\bar {\lambda}}_{v_i,\Gamma}(u)-1+\alpha t\pmod{I_p(G_1)I_p(H_1)}.$$
Therefore,
$${\det}_{G_1}=\Reg_{\Gamma}^{\Phi}t^n+b_{n+1}t^{n+1}+\dots.$$
This together with (\ref{e:rfirstterm}) and (\ref{e:aandb}) implies the theorem.

\end{proof}
\end{subsection}

\begin{subsection}{A twisted version}\label{sec:eqversion}
\begin{theorem}\label{thm:stark2}
Let notations be as in {\em{Corollary \ref{c:r}}}.
Then for every admissible $H$,
and every $\chi\in
{\hat \Gamma}$, the twisted Stickelberger element
$\theta_{\chi}\in I_{{\O},H}^n$. Furthermore,
there is a unique
$\epsilon\in \Z_{(p)}\Lambda_{S,T}^n$ such that for every admissible $H$ and
every $\chi \in {\hat {\Gamma}}$,
\begin{equation}\label{e:stark2}
\mathcal{R}_{\eta,\chi,H}=[\theta_{\chi}]_{(n,H)}.
\end{equation}
\end{theorem}
\begin{proof}
We only need to use Theorem \ref{thm:stark1} and then apply the $\chi$-twisted map to
the Stickelberger element and the refined regulator. The uniqueness is a consequence of
Lemma \ref{le:uniqueintegral}.
\end{proof}

\end{subsection}
\end{section}

\begin{section}{The decomposition of the refined regulator}\label{sec:decomporeg}

\begin{subsection}{The canonical pairing}
\label{sec:pairing}
In this chapter, we show that ${\det}_{\H}$, the regulator of Gross,
is decomposed into a product of some kind of irreducible factors.
To do so, it is
helpful to consider the dual version of the homomorphism
$\lambda_{H'}$, which can be described as a pairing.
% For simplicity, through $\delta_{H'}$, we will identify $H'$ with $I(H')/I(H')^2$.
\begin{definition}\label{def:pairing}
Let
$$Y'(K')=\Hom_{\Z} (Y(K'),\Z),$$
$$X'(K')=\Hom_{\Z} (X(K'),\Z),$$
and let
$$<\cdot , \cdot>_{H'}: U(K')\times
Y'(K')\longrightarrow I(H')/I(H')^2$$ be the pairing defined by
\begin{equation}
<u,\phi>_{H'}=\sum_{w\in S(K')} \phi(w) \cdot \lambda_w
(u),\;\text{for all }\;\; u\in U(K'),\;\phi\in Y'(K').
\label{eq:pairing}
\end{equation}
This pairing factors through a unique pairing on $U(K')\times
X'(K')$, which, by the abuse of notations, will also be denoted as
$<\cdot , \cdot>_{H'}$.
\end{definition}
\par
Following from (\ref{eq:compalambda}), for $H'\subset H''$, we have
\begin{equation}
<u,\phi\mid_{Y'(K'')}>_{H''}=<u,\phi>_{H'},\;\mbox{for all }\;\; u\in U(K'')
,\;\phi\in Y'(K').
\label{eq:compapairing}
\end{equation}

Directly from the definition, we have
\begin{equation}
<{}^{\sigma} u,{}^{\sigma} \phi>_{H'}=<u,\phi>_{H'},\;\mbox{for all}\; \;
\sigma\in G, u\in U(K'),\; \phi\in Y'(K'),
\label{eq:sigmapairing}
\end{equation}
where
$${}^{\sigma}\phi(w)=\phi(\sigma^{-1}(w)).$$

If $A=\{ a_i\}_{i=1,...,r_{K'}}$, $B=\{
b_j\}_{j=1,...,r_{K'}}$ are $\Z$-bases of $U(K')$ and $X'(K')$,
then as usual the associated discriminant of the pairing $<\cdot,\cdot>_{H'}$
is defined as
$$
\begin{array}{cl}
{} & \det(<a_i,b_j>)_{i,j=1,...,r_{K'}}\\
= & \sum_{\pi\in S_{r_{K'}}} \sign(\pi) <a_1,b_{\pi
(1)}>\cdot\dots\cdot <a_{r_{K'}},b_{\pi(r_{K'})}>,
\end{array}
$$
which is considered as an element of $I(H')^{r_{K'}}/I(H')^{r_{K'}+1}$.
This discriminant is independent of the choice of the bases up to $\pm 1$.
The following Lemma is obvious.
\begin{lemma}\label{l:discdet}
Let $w_1',...,w_{r_{K'}}'\in S(K')$ and the basis
$u_1,...,u_{r_{K'}}\in U(K')$ be as
in the definition of ${\det}_{H'}$. Choose the basis
$\phi_1,...,\phi_{r_{K'}}\in X(K')'$
such that
\newline
$\phi_i(a_1w_1'+\dots+a_iw_i'+\dots+a_{r_{K'}}w_{r_{K'}}')=
a_i$. Then the associated discriminant of $<\cdot,\cdot>_{H'}$
equals  $[{\det}_{H'}]_{(n)}$.
\end{lemma}

If $L/K$ is the constant $\Z_p$-extension and $\sigma\in H$ is the
Frobenius, then
\newline
$\Val_{\sigma,r}([{\det}_{H}]_{(r)})$ is just the
classical regulator of $U$. Therefore, the following lemma holds.
% which is the
%adjoin of $K$ with $\F_{q^{p^{\infty}}}$. Then the Artin map $\Art$
%factors through the degree map $\Sigma_K$, and we denote
%\begin{equation}
%\Art: \A_K^*  \stackrel{\Sigma_K}{\longrightarrow} \Z
%\stackrel{\sim}\longrightarrow H_0\;\subset H. \label{eq:degk}
%\end{equation}
\begin{lemma}\label{le:dual}
Suppose that $L/K$ is the constant $\Z_p$-extension. Then the
pairing
$$U\times X'(K):\stackrel{<\cdot , \cdot >_H}{\longrightarrow}
I(H)/I(H)^2\stackrel{\delta_H^{-1}}{\longrightarrow} H=\Z_p$$
has all its values in $\Z\subset H$. Furthermore, this pairing
induces  a perfect pairing
on $ \Q U \times \Q  X'(K)$.
\end{lemma}

The $\chi$-eigenspace of $F_p X'(K)$ will be denoted as
$X'_{\chi}$. We have
$$F_p X'(K)=\bigoplus_{\chi\in {\hat {\Gamma}}}
X'_{\chi}.$$
We denote by $<\cdot,\cdot>_{H,p}$ the induced pairing
on $F_p U\times F_p X'(K)$.
Note that by Equation (\ref{eq:sigmapairing}) and
Lemma \ref{le:dual}, we have a duality between
the $F_p[\Gamma]$-modules $U_{\chi}$ and
$X'_{\chi^{-1}}$, where $U_{\chi}$ is the
$\chi$-eigenspace of
$F_p\otimes U$ with $\dim_{F_p}U_{\chi}=r_{\chi}$.

\begin{definition}\label{def:pairchi}
Let $<\cdot,\cdot>_{\chi}$ denote the restriction of
$<\cdot,\cdot>_{H,p}$ to $U_{\chi}\times
X'_{\chi^{-1}}$.
If
$A_{\chi}=\{c_1,...,c_{r_{\chi}}\}$ and
$B_{\chi^{-1}}=\{d_1,...,d_{r_{\chi}}\}$ are $F_p$-bases of
$U_{\chi}$ and $X'_{\chi^{-1}}$, then define
$$
\begin{array}{rrlll}
[{\det}_{\chi}]_{(r_{\chi})} & = &
\det(<c_i,d_j>)_{i,j=1,...,r_{\chi}} & {} & {}\\
{} & = & \sum_{\pi\in S_{r_{\chi}}} \sign(\pi) <c_1,d_{\pi
(1)}>\cdot
...\cdot <c_{r_{\chi}},d_{\pi(r_{\chi})}>, & {}\\
\end{array}
$$
which is viewed as an element in
$I_{F_p}(H)^{r_{\chi}}/I_{F_p}(H)^{r_{\chi}+1}$.
\end{definition}

This discriminant is independent of the choice of bases
up to elements of $F_p^*$.
\end{subsection}

\begin{subsection}{$(K/k,S)$-extensions }\label{sec:ks}
\par
The main result of this chapter is Proposition
\ref{prop:factorization}, which concerns the decomposition of the
refined regulator ${\det}_{H'}$. For this purpose, we need to
consider Galois extensions that might not be admissible.
\begin{definition}\label{def:gsext}
Let ${\tilde L}/k$ be a Galois extension.
Then ${\tilde L}/k$ is called a $(K/k,S)$-extension,
if $K\subset {\tilde L}$ and the associated field extension
${\tilde L}/K$ is pro-$p$, abelian and unramified outside $S(K)$.
If ${\tilde L}$ contains the given field $L$, then we say
that ${\tilde L}/k$  is a $(K/k,S)$-extension
of $L/k$ and ${\tilde H}:=\Gal({\tilde L}/K)$ is a $(K/k,S)$-extension
of $H=\Gal(L/K)$. We say that ${\tilde L}/k$  is a strict $(K/k,S)$-extension of
$L/k$ and ${\tilde H}$ is a strict $(K/k,S)$-extension of $H$, if $L/k$ is
the maximal admissible field extension contained in ${\tilde
L}/k$.
\end{definition}
\par
If $\tilde H$ is a strict $(K/k,S)$-extension of $H$, then the natural
quotient map ${\tilde H}\longrightarrow H$ factors through

$${\tilde H}\longrightarrow {\tilde H}/\sum_{\gamma\in\Gamma} (1-\gamma)
{\tilde H} \longrightarrow H,$$
where the second arrow has a finite kernel.
\begin{lemma}\label{le:gsext}
Suppose that $L'/K$ is a pro-$p$ abelian extension unramified
outside $S(K)$. Then there exists an extension ${\tilde L}/k$, which
contains $L'/k$ and is a $(K/k,S)$-extension of $L/k$ . If both
$H$ and $\Gal(L'/K)$ are finitely generated over $\Z_p$, then
${\tilde L}/k$ can be chosen such that $\Gal({\tilde L}/K)$ is
finite free over $\Z_p$ and its maximal admissible quotient is
also finite free over $\Z_p$.
\end{lemma}

\begin{proof}

\par
For the first statement, we let ${\tilde L}$  be the adjoin of $L/K$
with all the fields ${}^{\sigma}L'$, $\sigma\in \Gal({\bar
k}^{sep}/k)$.
\par
If the conditions of the second statement hold, then $\Gal({\tilde
L/K})$ is finitely generated over $\Z_p$. Lemma \ref{le:leop} (for
${\K}=K$ and ${\cS}=S(K)$) then implies that there is an abelian
extension $M/K$ which is unramified outside $S(K)$, with Galois
group finite free over $\Z_p$ and contains ${\tilde L}/K$. Let
$M'$ be the adjoin of all the fields ${}^{\sigma}M$,
$\sigma\in \Gal({\bar k}^{sep}/k)$. Then $M'/k$ is a
$(K/k,S)$-extension and $Gal(M'/K)$ is finite free over
$\Z_p$.
%Let $M_0$ be the fixed field of the torsion part of
%$Gal(M'/K)$. Then $M_0/k$ is also a $(K/k,S)$-extension of $L/k$.
%Now, $\Gal(M_0/K)$ is finite free over $\Z_p$ and dominates
%${\tilde L}/K$.
The maximal admissible quotient of $Gal(M'/K)$
might not be free over $\Z_p$. By Lemma \ref{le:maxadmissible}, we
can chose an admissible extension $L''/k$ such that the Galois group
$Gal(L''/K)$ is finite free over $\Z_p$ and contains the maximal admissible
quotient of $Gal(M'/K)$. We complete the proof by replacing ${\tilde L}$
by $M'L''$.
\end{proof}
\end{subsection}
\begin{subsection}{Universal pairings} \label{sec:unipa}
\par
\begin{definition}\label{def:universal}
Let $A$, $B$ and $C$ be $F_p$-vector spaces and let $\{a_i\}_i$,
$\{b_j\}_j$ be bases of $A$ and $B$. A pairing  of $F_p$-spaces
$$\Phi:A\times B\longrightarrow C$$
is said to be  universal, if the set $\{\Phi(a_i,b_j)\}_{i,j}$ is
linearly independent over $F_p$.
\end{definition}
\par
In other words, $\Phi$ is universal if and only if its image generates
a subspace which is canonically
isomorphic to $A\otimes_{F_p} B$.
\par
Suppose that ${\tilde L}/k$ is a $(K/k,S)$-extension and $\tilde
H=\Gal({\tilde L}/K)$. Then $\tilde H$ is abelian. We can use
(\ref{eq:pairing}) and define the paring
$$<\cdot,\cdot>_{\tilde H}:U\times X'(K)\longrightarrow I({\tilde
H})/I({\tilde H})^2.$$
Then the pairing $<\cdot,\cdot>_{\tilde H}$ is
$\Gamma$-equivariant in the sense that
\begin{equation}
<{}^{\gamma} u,{}^{\gamma} \phi>_{\tilde H}={}^{\gamma }<u,\phi>_{\tilde
H},\;\mbox{for all}\;\; u\in U, \phi\in X'(K), \gamma\in\Gamma.
\label{eq:gammaequi}
\end{equation}
\par
\begin{definition}
\label{def:univhtilde} A $(K/k,S)$-extension ${\tilde L}/k$ is
universal if it satisfies the following.
\begin{enumerate}
\item The Galois group $\tilde H=\Gal({\tilde L}/K)$ is finite
free over $\Z_p$.
\item The induced pairing
$$F_p\cdot U \times F_p\cdot  X'(K)\longrightarrow F_p\cdot
I({\tilde H})/I({\tilde H})^2=I_{F_p}(\tilde H)/I_{F_p}(\tilde H)^2$$
is universal.
\end{enumerate}
\end{definition}

\begin{definition}\label{def:unrestricted}
An admissible extension $L/k$, as well as the Galois group
$H=\Gal(L/K)$, is said to be unrestricted, if the following
conditions are satisfied.
\begin{enumerate} \item $H\simeq\Z_p^d$, for some $d$.
\item There exists a universal strict $(K/k,S)$-extension of
$H$.
\item The constant $\Z_p$-extension $L_0/K$ is a sub-extension of $L/K$.
\end{enumerate}
\end{definition}

\begin{lemma}
\label{le:unrestricted}
The followings are true.
\begin{enumerate}
\item Suppose that ${\tilde L}'/k$ contains ${\tilde L}/k$ and they are
both $(K/k, S)$-extensions. If ${\tilde H}':=\Gal({\tilde L}'/K)$
is finite free over $\Z_p$ and ${\tilde L}/k$ is universal, then
${\tilde L}'/k$ is also universal.
\item Suppose that $L'/k$ contains $L/k$ and they are both admissible. If
$H':=\Gal(L'/K)$ is finite free over $\Z_p$ and $L/k$ is
unrestricted, then $L'/k$ is also unrestricted.
\item Every admissible extension whose Galois group is finitely generated over $\Z_p$
is contained in certain unrestricted admissible extension.
\end{enumerate}
\end{lemma}
\begin{proof}
\par
Part (1) follows directly from the definitions. For Part (2), we
just note that if ${\tilde L}/k$ is a universal strict
$(K/k,S)$-extension of $L/k$, then Part (1) implies that
${\tilde L}L'/K$ is a universal strict $(K/k,S)$-extension
of $L'/k$.
\par
To prove Part (3), we first note that by Lemma 3.3 of \cite{tn}, a
universal $(K/k,S)$-extension $M/k$ exists. We denote by $L'$ the
field obtained by adjoining $M$ with the given admissible
extension and the constant $\Z_p$-extension. By Lemma \ref{le:gsext},
there is a $(K/k,S)$-extension
${\tilde L}/k$ such that ${\tilde L}$ contains $L'$ and both
$\Gal({\tilde L}/K)$ and its maximal admissible quotient are
finite free over $\Z_p$. Since ${\tilde L}$ contains $M$ and
$M/k$ is universal, by Part (1), ${\tilde L}/k$ is also universal.
This completes the proof.
\end{proof}
\end{subsection}

\begin{subsection}{The decomposition}

\label{sec:factorization}

\par
In view of
(\ref{eq:sigmapairing}), for the admissible Galois group $H$, the pairing
$<\cdot,\cdot>_{H,p}$ is $\Gamma$-invariant  and hence can not be universal. However, we
are going to show that if $H$ is unrestricted, then its $\chi$-part,
$< \cdot , \cdot>_{\chi}$, is universal for every $\chi\in {\hat {\Gamma}}$.
\par
\begin{lemma}\label{le:chiuniversal}
Assume that $H$ is  unrestricted. The followings are true.
\begin{enumerate}
\item For each $\chi\in {\hat {\Gamma}}$, the pairing
$<\cdot,\cdot>_{\chi}$ is universal over $F_p$.
\item The linear sub-spaces
$<U_{\chi},X'_{\chi^{-1}}>_{\chi}$, $\chi\in {\hat {\Gamma}}$, of
$F_p\cdot I(H)/I(H)^2$ are linearly independent over $F_p$.
\end{enumerate}
\end{lemma}
\begin{proof}
\par
Let $\tilde H$ be a universal strict $(K/k,S)$-extension of $H$.
Let the ${\tilde H}^{(1)}\subset \tilde H$ be the $1$-eigenspace of $\Gamma$.
Then the natural quotient map ${\tilde
H}\longrightarrow H$ induces an isomorphism
\begin{equation}
  F_p \otimes_{\Z_p} {\tilde H}^{(1)}{\longrightarrow} F_p \otimes_{\Z_p} H.
\label{eq:h1}
\end{equation}
To simplify the notations, for the rest of the proof, for every
Galois group $\mathcal{H}$, through the isomorphism
$\delta_{\mathcal{H}}$ (see (\ref{eq:del})), we will identify
$I(\mathcal{H})/I(\mathcal{H})^2$ with $\mathcal{H}$ .

By the $\Gamma$-equivariant property of the pairing
$<\cdot,\cdot>_{\tilde H}$, the restriction of \newline $id_{F_p}
\otimes <\cdot,\cdot>_{\tilde H}$ to $U_{\chi}\times
X'_{\chi^{-1}}$ has values all in $F_p\otimes_{\Z_p} {\tilde
H}^{(1)}$.
\par
To say that $id_{F_p}\otimes <\cdot,\cdot>_{\tilde H}$ is
universal is the same as to say that the induced homomorphism
$F_p\otimes U\otimes X'\longrightarrow F_p\otimes_{\Z_p} {\tilde
H}$ is an injection. Since $\bigoplus_{\chi\in {\hat {\Gamma}}}
U_{\chi}\otimes X'_{\chi^{-1}}$ is a direct summand of $F_p\otimes
U\otimes X'$ and its image under the induced homomorphism is in
$F_p\otimes_{\Z_p} {\tilde H}^{(1)}$, the lemma is proved by
taking the isomorphism (\ref{eq:h1}).
\end{proof}

\end{subsection}

\begin{subsection}{The associated homogeneous polynomials}
\label{sec:po}
Let $\Ver$ denote the transfer homomorphism
\begin{equation}\label{e:transfer}
\begin{array}{rrcl}
\Ver: & G & \longrightarrow & H\\
{} & g & \mapsto & |\Gamma|g.\\
\end{array}
\end{equation}
For simplicity, we will also let $\Ver$ denote the restriction of it
to a subgroup $H'$ as well as the induced maps on augmentation quotients.

\par
\begin{definition}\label{def:po}
Suppose that $H\simeq\Z_p^d$ and $\E$ is a $\Z_p$ basis. Let $R=\O_p$ and
recall the notations in {\em{(\ref{eq:mapd})}}, {\em{(\ref{eq:mapdf})}}.
For each $H'$ define the homogeneous polynomial
$$f_{H'}:=d_{{\E},F_p}([\Ver({\det}_{H'})]_{r_{K'}})
)\in F_p[s_1,...,s_d].$$ Also, for each $\chi\in {\hat
{\Gamma}}$, define the homogeneous polynomial
$$
f_{\chi}:=f_{{\E},\chi}:=d_{{\E},F_p}([\Ver({\det}_{\chi})]_{r_{\chi}}) \in
F_p[s_1,...,s_d].$$
\end{definition}
\par
The polynomial $f_{H'}$ is independent of the choice of the basis
of $H$ up to elements of $\Z_p^*$. The polynomial $f_{\chi}$ is
uniquely defined up to
elements of $F_p^*$.
%\par
%Recall the degree map $\A_K^* \stackrel{\Sigma_K}{\longrightarrow}
%\Z$ described in (\ref{eq:degk}). If $H$ is unrestricted, then the
%degree map factors through the Artin map $\A_K^*\longrightarrow
%H$. Let $\pi:H\longrightarrow \Z$ be the induced degree map. Also,
%to simplify the notations, we will identify $I(H)/I(H)^2$ with $H$,
%through the isomorphism $\delta_H$ (see (\ref{eq:del})).
%\begin{definition}
%\label{def:rational}
\begin{definition}\label{d:rational}
Let $H$ be an unrestricted admissible group. Then a $\Z_p$-basis ${\E}$ is called rational
if the following conditions are satisfied.
\begin{enumerate}
\item
\begin{equation}
<U,X'(K)>_{H} \subset \sum_{\sigma\in {\E}} \Z\cdot \delta_H(\sigma).
\label{eq:rational}
\end{equation}
\item
If $L_0/K$ is the constant $\Z_p$-extension with $\Gal(L_0/K)=H_0=\Z_p$
and $\pi:H\longrightarrow H_0$ is the natural projection, then
$$\pi(\sigma)\in\Q\cap\Z_p \subset\Z_p,\;\mbox{for all}\; \sigma\in {\E}.$$
\end{enumerate}
\end{definition}
\begin{lemma}\label{l:rational}
An unrestricted admissible Galois group always has a rational basis.
\end{lemma}
\begin{proof}
Suppose that $H\simeq \Z_p^d$ is unrestricted. Let $A$ be a $\Z$-basis of $U$ and
let $B$ be a $\Z$-basis of $X'(K)$. We can find a $d$-dimensional
$\Q$-vector space $V\subset \Q_p\otimes H$ such that $\{ <a,b>_H|\; a\in A, b\in B\}
\subset \delta_H(V)$.
Then $M:=V\cap H$ is a free $\Z$-module of rank $d$. Let $\E$ be a basis of $M$.
Then obviously, the inclusion (\ref{eq:rational}) holds.
Lemma \ref{le:dual} implies that $\E$ is rational.
\end{proof}

\begin{proposition}\label{prop:factorization}
Suppose that $H\simeq \Z_p^d$ is unrestricted.
and $\E$ is a rational basis.
Then for each $H'$, the homogeneous
polynomial $f_{H'}$ is in $\Q [s_1,...,s_d]$
and for every $\chi$, the polynomial $f_{\chi}$ is
absolutely irreducible and can be chosen in
$ F[s_1,...,s_d]$ such that for some number $c\in F^*$
\begin{equation}\label{eq:fac}
f_{H'}=c\cdot \prod_{\chi\in {\hat {\Gamma'}}} f_{\chi}.
\end{equation} Furthermore, the polynomials
$\{ f_{\chi}\; | \; \chi\in {\hat {\Gamma}}\} $
are algebraically independent over $F_p$.

\end{proposition}
\begin{proof}
The inclusion (\ref{eq:rational}) implies that $f_H\in \Q[s_1,...,s_d]$.
Similarly, in $\sum_{u\in U}
F\cdot u$, we can find a
$\Z_p$-basis of $U_{\chi}$ and in
$\sum_{x\in X'(K)} F\cdot x$ we can find a $\Z_p$-basis of $X'_{\chi^{-1}}$.
Then we have $f_{\chi}\in F[s_1,...,s_d]$.

Since $H$ is fixed by $\Gamma$ and $<\cdot,\cdot>_H$ is
$\Gamma$-invariant, for $\chi'\not=\chi^{-1}$ the restriction of
$ <\cdot,\cdot>_{H,p}$ to the set $U_{\chi}\times
X'_{(\chi')^{-1}}$ is the trivial pairing. This shows the
existence of Equation (\ref{eq:fac}) for $H'=H$. In general,
we use the compatibility equality (\ref{eq:compapairing}) and view
$<\cdot,\cdot>_{H',p}$ as a part of $<\cdot,\cdot>_{H,p}$.
\par
For each integer $m$, consider the determinant of the $m$ by $m$
matrix $(t_{ij})$, where $t_{ij},\; i=1,...,m,\; j=1,..., m$ are
independent variables over a field. It is well known that the
determinant of this matrix is absolutely irreducible (see
\cite{v}). This fact and Lemma \ref{le:chiuniversal} imply the
irreducibility and the algebraic independence of $f_{\chi}$'s.
\end{proof}

\end{subsection}

\begin{subsection}{The explicit expressions}
\label{sec:explicit}
Let $\eta$ and $w_1,...,w_n$ be as in Definition \ref{d:fix} and
let $w^{(i)}\in X'(K)$ be such that $w^{(i)}(\sum_{w\in S(K)}\alpha_ww)=
\alpha_{w_i}$, for every
$ \sum_{w\in S(K)}\alpha_ww\in X'(K)$.
Then $\{w^{(1)},..., w^{(n)} \}$ generate
a free $\Z[\Gamma]$-module of rank $n$. Assume that $\chi$ is a character
of $\Gamma$ such that $r_{\chi}=n$. Then
$$\{\frac{1}{|\Gamma|} \sum_{\gamma\in\Gamma} \chi(\gamma)\cdot
^{\gamma}w^{(i)}|\;\;i=1,...,n\}$$ is a basis of the $F$-vector
space $X'(K)_{\chi^{-1}}$. Also, if $M$ is a $\Q[\Gamma]$-free sub-space of
$\Q\cdot U$ of rank $n$ and is generated over
$\Q[\Gamma]$ by $\{\epsilon_1,...,\epsilon_n\}$, then
$$\{\sum_{\gamma\in\Gamma} \chi(\gamma^{-1})\cdot ^{\gamma}\epsilon_i\;|\;
i=1,...,n\}.$$ is a basis of the $F_p$-vector space $U_{\chi}$.
Using  Definition \ref{def:twistregulator} and equations (\ref{e:iotaH})and
(\ref{eq:sigmapairing}), we obtain
\begin{equation}\label{le:explicit}
\begin{array}{rcl}
[\Ver({\det})_{\chi}]_{r(\chi)} &=& \Ver \Bigl({\det}
\Bigl(\sum_{\gamma\in\Gamma}\chi(\gamma)\cdot
\lambda_{\gamma(w_j)}(\epsilon_i)\Bigr)_{\stackrel{i=1,...,n}{j=1,...,n}}
\Bigr)\\
{} &=& \Ver(\mathcal{R}_{\eta,\chi,H}^{\triangleright}
(\epsilon_1\wedge ...\wedge \epsilon_n))
\end{array}
\end{equation}
%Thus, we have proven the following lemma.
%\begin{lemma}\label{le:explicit}
%Suppose that $\chi$ is a character of $\Gamma$ such that
%$r(\chi)=n=sp(K/k)$ and $v_1,...,v_n\in S$ split completely in
%$K$. For each $i$ fix a place $w_i$ of $K$ such that $w_i|v_i$ and
%put $\eta=w^*_1\wedge ... \wedge w^*_n$. Then there is an
%$\epsilon\in\Q \Lambda^n U$ such that
%$${\rm det}_{\chi}=t(\mathcal{R}eg_{\eta,\chi,H}(\epsilon)).$$
%\end{lemma}
\end{subsection}
\end{section}

\begin{section}{The proof}\label{sec:proof}

\begin{subsection}{The product formula}\label{sec:prodform}
\par
As before, $H'$ is a subgroup of $G$ such that $H\subset H'$.
Since $H'$ is a closed subgroup, each character $\psi\in {\hat
H'}$ can be extended to a character on $G$. Let ${\hat
{\Gamma'}}_{\psi}$ denote the set of all such extensions of $\psi$.
Recall the definitions in Section \ref{sec:theta} of the modified
$L$-function $L_{S(K'),T(K')}(\psi,s)$ and the Stickelberger
element $\theta_{H'}$. By Class Field Theory, we have the
following product formula.
\begin{equation}
L_{S(K'),T(K')}(\psi,s)=\prod_{\phi\in {\hat {\Gamma'}}_{\psi}}
L_{S,T}(\phi,s). \label{eq:prodl}
\end{equation}
\begin{proposition}\label{prop:thetafac}
Suppose that $H'$ is a subgroup of G such that $H\subset H'$. Then in
the group ring $\O[G]$, we have
\begin{equation} \theta_{H'}=\prod_{\chi\in {\hat {\Gamma'}}} \theta_{\chi}.
\label{eq:thetaprod}
\end{equation}
\end{proposition}
\begin{proof}
\par
Apply every $\psi$ to both sides of (\ref{eq:thetaprod}), then use the product
formula (\ref{eq:prodl}).
\end{proof}
\par
\begin{definition}\label{def:thetapo}
Assume that $H\simeq\Z_p^d$ and ${\E}=\{\sigma_1,...,\sigma_d\}$
is a basis of $H$ over $\Z_p$.
Let $R=\O_p$ and
recall the notations in {\em{(\ref{eq:mapd})}}, {\em{(\ref{eq:mapdf})}}.
For each $H'$
define
$$\xi_{H'}=d_{{\E},F_p}([\Ver(\theta_{H'})]_{(r_{K'})}).$$
For each $\chi\in {\hat {\Gamma}}$, let $n_{\chi}$ be such that
$\Ver(\theta_{\chi})\in I_{\O_p}(H)^{n_{\chi}}\setminus
I_{\O_p}(H)^{n_{\chi}+1}$ and define
$$\xi_{\chi}=d_{{\E},F_p}([\Ver(\theta_{\chi})]_{(n_{\chi})}).$$
\end{definition}

\begin{proposition}\label{prop:thetafac2}
Suppose that $H\simeq \Z_p^d$ is unrestricted
and $\E$ is a rational basis.
Then for each $H'$, the homogeneous
polynomial
\begin{equation}\label{e:gross2}
\xi_{H'}=h_{K',S(K'),T(K')}f_{H'}.
\end{equation}
In particular, $\xi_{H'}$ is in $\Q [s_1,...,s_d]$ and is nonzero.
Furthermore, we have
\begin{equation}\label{eq:thetafac}
\xi_{H'}=\prod_{\chi\in {\hat {\Gamma'}}} \xi_{\chi}.
\end{equation}
\end{proposition}
\begin{proof}
It is a consequence of Theorem \ref{th:gross} and Proposition
\ref{prop:thetafac}
\end{proof}
\end{subsection}

\begin{subsection}{The end of the proof}
\label{sec:end}
\par
In this section, we complete the proof of the main theorem.

\vskip5pt
\begin{proof}
We only need to prove the theorem for the case where $H=\H$. Since $\H$ is the projective limit
of unrestricted admissible groups,
we will first consider the case where $H$ is unrestricted and $\E$ is a rational basis.
For $\chi_1,\chi_2\in
\hat \Gamma$, denote $\chi_1\sim\chi_2$ if they generate the same
cyclic subgroup of $\hat\Gamma$.
\vskip7pt \noindent {\bf {Step 1}}:\begin{em} For each $\chi\in\hat\Gamma$,
there is a character $\chi'\sim\chi$ and a constant $c(\chi,\chi')\in F^*$
such that
\begin{equation}
\xi_{\chi}=c(\chi,\chi')\cdot f_{\chi'}.
\label{eq:chichi'}
\end{equation}
\end{em}
\vskip7pt
\par
To show this, we recall equations (\ref{eq:fac}) ,(\ref{e:gross2})
and (\ref{eq:thetafac}) and let $H'$ run through all the subgroup of $G$
containing $H$. Consequently, there is a $c_{\chi}\in F_p^*$ such that
$$\prod_{\chi'\sim\chi}\xi_{\chi'}=c_{\chi}\cdot \prod_{\chi'\sim\chi}
f_{\chi'}.$$

Under the natural action of $Gal({\bar {\Q}}/{\Q})$ on ${\hat
{\Gamma}}$, the set $\{\chi'| \; \chi'\sim\chi\}$ form an orbit.
For each $\tau\in Gal({\bar {\Q}}/{\Q})$, we have $\theta_{{}^{\tau}\chi}
={}^{\tau}\theta_{\chi}\in F[G]$ and hence
$$n_{\chi'}=\deg(\xi_{\chi'})=\deg(\xi_{\chi})=\deg (f_{\chi})=r_{\chi}, \;\text{if} \;\chi'\sim\chi.$$

Since all
the $f_{\chi'}$ are absolutely irreducible and algebraically
independent, there is a $\chi'$ and a number $c(\chi,\chi')$ such that (\ref{eq:chichi'})
holds.

We need to show $c(\chi,\chi')\in F^*$. The problem is that we don't know if the polynomial
$\xi_{\chi}$ has its coefficients in $F$, although by Proposition \ref{prop:factorization}
this holds for the polynomial $f_{\chi'}$. To overcome this problem, we apply the natural
projection
$$\pi:H\longrightarrow H_0=\Gal(L_0/K),$$
where, as before, $L_0$ is the constant
$\Z_p$-extension of $K$. Let $\sigma$ be the Frobenius
in $H_0$ and let $t=\sigma-1$. Then $\pi(\E)\subset \Q\sigma$ and
$\pi_*(f_{\chi'})=b_{r_{\chi'}}t^{r_{\chi}}$ with $b_{r_{\chi'}}\in F$ for each $\chi'$.
Since $\pi_*(f_H)=b_Ht^{r}$, where $b_H$ is a nonzero multiple of the classical regulator
of $U$, by the product formula (\ref{eq:fac}), we have $b_{r_{\chi'}}\in F^*$.

Also, Lemma
\ref{l:Thetatheta}(2) and the
functorial property of the Stickelberger elements imply that
if $\pi_*(\Ver(\theta_{\chi}))=a_{r_{\chi}}t^{r_{\chi}}+\dots\in\F_p[[t]]$,
then $a_{r_{\chi}}$ is a  nonzero rational multiple of the $r_{\chi}$th derivative
$L_{S,T}(\chi,0)^{(r_{\chi})}\in F^*$.
Therefore
$$c(\chi,\chi')=a_{r_{\chi}}/b_{r_{\chi}}\in F^*.$$

\vskip7pt
\noindent
{\bf {Step 2}}: $\chi'=\chi$.
\vskip7pt
\par
To show this, we need to apply a work of Hayes. According to
\cite{h}, our Theorem \ref{thm:stark1} and its consequence
Theorem \ref{thm:stark2} are true in the case where $n=1$.
Let
$S^{(1)}=S\setminus\{v_2,...,v_n\}$. Since $v_1,...,v_n\in S$ split
completely in $K$, the extension $K/k$ is unramified outside $S^{(1)}$. Also, $\#
S^{(1)}\geq 2$, since $n\leq \# S-1$. If an abelian extension $L^{(1)}/K$ with Galois group
$H^{(1)}$ is admissible with respect to the setting $(K/k,S^{(1)})$, then it is
also admissible with respect to $(K/k,S)$. We assume that $H^{(1)}$ is
unrestricted with respect to $(K/k,S')$ and (after certain
extension, if necessary) the given $L$ contains $L^{(1)}$. We will use
$\theta^{(1)}$, $\xi^{(1)}$, $\det^{(1)}$, $f^{(1)}$, $\xi^{(1)}$
and so on to denote the objects
derived from $H^{(1)}$ and $S^{(1)}$. Let $Q:H\longrightarrow H^{(1)}$ be the
natural quotient map. Then for every $\psi\in {\hat {\Gamma}}$
$$Q_*(\theta_{\psi})=(1-[w_2])\dots (1-[w_n])\cdot \theta^{(1)}_{\psi},$$
where $[w_i]\in H'$ is the Frobenius element at $w_i$. Also, there
is a $c'(\psi)\in F^*$ such that
$$Q_*([{\det}_{\psi}]_{(r_{\psi})})=c'(\psi)\cdot (1-[w_2])\dots (1-[w_n])\cdot
[{\det}^{(1)}_{\psi}]_{(r_{\psi}-n+1)}.$$ Hayes' result together with Equation
(\ref{le:explicit}) implies that
\begin{equation}
Q_*(f_{\chi})=c''(\chi)\cdot Q_*(\xi_{\chi}),
\label{eq:chichi}
\end{equation}
for some $c''(\chi)\in F^*$.
Equation (\ref{eq:chichi'}) and (\ref{eq:chichi}) imply that
$f^{(1)}_{\chi'}$ and $f^{(1)}_{\chi}$ are proportional to each other. This can not
happen unless $\chi'=\chi$, since $H^{(1)}$ is unrestricted and all the
$f^{(1)}_{\psi}$, $\psi\in {\hat {\Gamma}}$, are algebraically independent.
\vskip7pt
\noindent
{\bf {Step 3}}: There is an $\epsilon$ in $\Q\Lambda ^n U $ such that for every
$ \chi\in {\hat {\Gamma}}$,
\begin{equation}
[\Ver(\theta_{\chi})]_{(n)}=\Ver( \mathcal{R}_{\eta,\chi,H}^{\triangleright}(\epsilon)).
\label{eq:c'cc}
\end{equation}
\noindent
\vskip7pt

By (\ref{le:explicit}), there is an
$\epsilon'$ in $\Q\Lambda ^{n} U$ such that for each $\chi\in\hat\Gamma$
there is a $c_{\chi}\in F^*$ such that
$$[\Ver(\theta_{\chi})]_{(n)}=c_{\chi}\cdot
\Ver(\mathcal{R}_{\eta,\chi,H}^{\triangleright}(\epsilon'))\in
I_{F_p}(H)^n/I_{F_p}(H)^{n+1}.$$
Here we have
\begin{equation}\label{e:lambdast}
c_{\chi}=0,\; \text{for a character}\; \chi \; \text{such that}\; r_{\chi}>n,
\end{equation}
and also, for every
$\tau\in \Gal(\bar\Q/\Q)$,
$$c_{{}^{\tau}\chi}={}^{\tau}c_{\chi}.$$
Then there is an element $\alpha\in\Q[\Gamma]$ such that $\chi(\alpha)=c_{\chi}$
for every $\chi\in\hat\Gamma$. Let $\epsilon=\alpha\cdot\epsilon'$. Then
Equation (\ref{eq:c'cc}) holds for every $\chi$.
\vskip7pt
\noindent
{\bf {Step 4}}: $\theta_G\in I_H^n$.
\vskip7pt
\par

It is easy to see that
\begin{equation}\label{e:ver}
\begin{array}{crcl}
\Ver: & I_{F_p}(H)^m/I_{F_p}(H)^{m+1} & \longrightarrow & I_{F_p}(H)^m/I_{F_p}(H)^{m+1}\\
{} & x & \mapsto & |\Gamma|^m x\\
\end{array}
\end{equation}
is injective.

As in definition \ref{def:retheta}, for each $\gamma\in\Gamma$,
let $\theta_{\gamma}$ be the $\gamma$-part of $\theta_G$. Then
$$\Ver(\theta_{\chi})=\sum_{\gamma\in\Gamma}\chi(\gamma)\cdot\Ver(\theta_{\gamma}),$$
and Equation (\ref{eq:c'cc}) implies that each $\Ver(\theta_{\gamma})\in I_{F_p}(H)^n$.
If $\theta_{\gamma}\in I_{F_p,H}^m\setminus  I_{F_p,H}^{m+1}$,
then $\pounds^{-1}([\theta_{\gamma}]_{(m,H)})=\gamma\otimes \theta_{[\gamma],m}\not=0$, for some
$\theta_{[\gamma],m}\in I_{F_p}(H)^m/I_{F_p}(H)^{m+1}$ and hence
\newline
$[\Ver(\theta_{\gamma})]_{(m,H)}=|\Gamma|^m\theta_{[\gamma],m}\not=0$.
But we also have $[\Ver(\theta_{\gamma})]_{(m,H)}=[\Ver(\theta_{\gamma})]_{(m)}$,
which is zero unless $m\geq n$. This shows every $\theta_{\gamma}$ is in
$I_{F_p,H}^n$ and so is $\theta_G$. Lemma \ref{le:basechange} and Lemma \ref{le:torfree}
imply that $\theta_G\in I_H^n$.
\vskip7pt
\noindent
{\bf {Step 5}}: $[\theta_G]_{(n,H)}=\pounds(\mathcal{R}_{\eta,H}^{\triangleright}(\epsilon))$.
\vskip7pt
\par

Put
$$\mathcal{R}_{\eta,H}^{\triangleright}(\epsilon)
=\sum_{\gamma\in \Gamma} \gamma\otimes \mathcal{R}_{\gamma},\;\;
\mathcal{R}_{\gamma}\in I_{F_p}(H)^n/I_{F_p}(H)^{n+1},$$
and
$$\pounds^{-1}([\theta_G]_{(n,H)})=\sum_{\gamma\in \Gamma} \gamma\otimes \theta_{[\gamma]},\;\;
\theta_{[\gamma]}\in I_{F_p}(H)^n/I_{F_p}(H)^{n+1}.$$
Then
$$\Ver(\mathcal{R}_{\eta,\chi,H}^{\triangleright}(\epsilon))
=\sum_{\gamma\in \Gamma}|\Gamma|^n\chi( \gamma)\mathcal{R}_{\gamma}\in I_{F_p}(H)^n/I_{F_p}(H)^{n+1},$$
and
$$[\Ver\theta_{\chi}]_{(n,H)}=\sum_{\gamma\in \Gamma} |\Gamma|^n\chi(\gamma) \theta_{[\gamma]}
\in I_{F_p}(H)^n/I_{F_p}(H)^{n+1}.$$
By applying the inverse Fourier transform to (\ref{eq:c'cc}), we
deduce that
$\theta_{[\gamma]}=\mathcal{R}_{\gamma}\in I_{F_p}(H)^n$ for every
$\gamma\in\Gamma$, and hence
\begin{equation}\label{e:cong}
[\theta_G]_{(n,H)}=\pounds(\mathcal{R}_{\eta,H}^{\triangleright}(\epsilon)).
\end{equation}
We have proved the above equality in the case where the coefficient ring is $F_p$,
but Lemma \ref{le:basechange} and Lemma \ref{le:torfree} say that the same equality holds
in the case where the
coefficient ring  is $\Z$.
\vskip7pt
\noindent
{\bf {Step 6}}: the case $H=\H$.
\noindent
\vskip7pt

Lemma \ref{le:uniqueintegral} and Lemma \ref{le:unrestricted} say that if $H$ is unrestricted
and is large enough then the exterior product $\epsilon$ in (\ref{e:cong}) is unique. Taking
projective limit, we see from the functorial properties that this $\epsilon$ also
satisfies (\ref{e:cong}) for the case where $H=\H$. Corollary \ref{c:lambda0H} says that
$\epsilon\in \Z_{(p)}\Lambda_0^n U$. From (\ref{e:lambdast}), we see that
$\epsilon\in
\Z_{(p)}\Lambda^n_{S,T}$.

\end{proof}

\end{subsection}

\end{section}

\end{document}